\documentclass[a4,11pt]{article}
\usepackage{amsmath}
\usepackage{graphics}
\usepackage{latexsym}
\usepackage{amsfonts}
\numberwithin{equation}{section}
\renewcommand{\baselinestretch}{1.40}
\title{Optimal Multi-Modes Switching Problem in Infinite Horizon}
\author{Brahim EL ASRI \thanks{Universit\'e du
Maine, D\'epartement de Math\'ematiques, Equipe Statistique et
Processus, Avenue Olivier Messiaen, 72085 Le Mans, Cedex 9, France.
e-mail: brahim.elasri@univ-lemans.fr}}
\begin{document}
\date{\today}
\maketitle
\newtheorem{theo}{Theorem}
\newtheorem{problem}{Problem}
\newtheorem{pro}{Proposition}
\newtheorem{cor}{Corollary}
\newtheorem{axiom}{Definition}
\newtheorem{rem}{Remark}
\newtheorem{lem}{Lemma}
\newtheorem{app}{Appendix}
\newcommand{\brm}{\begin{rem}}
\newcommand{\erm}{\end{rem}}
\newcommand{\beth}{\begin{theo}}
\newcommand{\eeth}{\end{theo}}
\newcommand{\bl}{\begin{lem}}
\newcommand{\el}{\end{lem}}
\newcommand{\bp}{\begin{pro}}
\newcommand{\ep}{\end{pro}}
\newcommand{\bcor}{\begin{cor}}
\newcommand{\ecor}{\end{cor}}
\newcommand{\be}{\begin{equation}}
\newcommand{\ee}{\end{equation}}
\newcommand{\beq}{\begin{eqnarray*}}
\newcommand{\eeq}{\end{eqnarray*}}
\newcommand{\beqa}{\begin{eqnarray}}
\newcommand{\eeqa}{\end{eqnarray}}
\newcommand{\dg}{\displaystyle \delta}
\newcommand{\cm}{\cal M}
\newcommand{\cF}{{\cal F}}
\newcommand{\cR}{{\cal R}}
\newcommand{\bF}{{\bf F}}
\newcommand{\tg}{\displaystyle \theta}
\newcommand{\w}{\displaystyle \omega}
\newcommand{\W}{\displaystyle \Omega}
\newcommand{\vp}{\displaystyle \varphi}
\newcommand{\ig}[2]{\displaystyle \int_{#1}^{#2}}
\newcommand{\integ}[2]{\displaystyle \int_{#1}^{#2}}
\newcommand{\produit}[2]{\displaystyle \prod_{#1}^{#2}}
\newcommand{\somme}[2]{\displaystyle \sum_{#1}^{#2}}
\newlength{\inter}
\setlength{\inter}{\baselineskip} \setlength{\baselineskip}{7mm}
\newcommand{\no}{\noindent}
\newcommand{\rw}{\rightarrow}
\def \ind{1\!\!1}
\def \R{I\!\!R}
\def \cadlag {{c\`adl\`ag}~}
\def \esssup {\mbox{ess sup}}
\begin{abstract}
This paper studies the problem of the deterministic version of the
Verification Theorem for the optimal $m$-states switching in
infinite horizon under Markovian framework with arbitrary switching
cost functions
. The problem is formulated as an extended impulse control problem
and solved by means of probabilistic tools such as the Snell envelop
of processes and reflected backward stochastic differential
equations. A viscosity solutions approach is employed to carry out a
fine analysis on the associated system of $m$ variational
inequalities with inter-connected obstacles. We show that the vector
of value functions of the optimal problem is the unique viscosity
solution to
the system. 
 This problem is in relation with the valuation of firms in
a financial market.
\end{abstract}

\no{\bf AMS Classification subjects}: 60G40 ; 62P20 ; 91B99 ; 91B28
; 35B37 ; 49L25.
\medskip

\no {$\bf Keywords$}: Real options; Backward stochastic differential
equations; Snell envelope; Stopping times ; Switching; Viscosity
solution of PDEs; Variational inequalities.

\section {Introduction}First let us deal with an example in order to introduce the problem we consider
in this paper:

Assume we have a power station/plant which produces electricity and
which has several modes of production, e.g., the lower, the middle
and the intensive modes. The price of electricity in the market,
given by an adapted stochastic process $(X_t)_{t \geq0}$, fluctuates
in reaction to many factors such as demand level, weather
conditions, unexpected outages and so on. On the other hand,
electricity is non-storable, once produced it should be almost
immediately consumed. Therefore, as a consequence, the station
produces electricity in its instantaneous most profitable mode known
that when the plant is in mode $i\in {\cal I}$, the yield per unit
time $dt$ is given by means of $\psi_i(X_t)dt$ and, on the other
hand, switching the plant from the mode $i$ to the mode $j$ is not
free and generates expenditures given by $g_{ij}(X_t)$ and possibly
by other factors in the energy market. The switching from one regime
to another one is realized sequentially at random times which are
part of the decisions. So the manager of the power plant faces two
main issues:

$(i)$ when should she decide to switch the production from its
current mode to another one?

$(ii)$ to which mode the production has to be switched when the
decision of switching is made?
\medskip

\noindent In other words she faces the issue of finding the optimal
strategy of management of the plant. This issue is in relation with
the price of the power plant in the energy market.
\medskip

Optimal switching problems for stochastic systems were studied by
several authors (see e.g. \cite{[BE], [BO1],[BS], [CL], [DP], [DH],
[DHP], [DZ2], [EH], [HJ], [LPZ], [TY], dz} and the references
therein). The motivations are mainly related to decision making in
the economic sphere. Several variants of the problem we deal with
here, including finite and infinite horizons, have been considered
during the recent years. In order to tackle those problems, authors
use mainly two approaches. Either a probabilistic one \cite{[DH],
[DHP], [HJ]} or an approach which uses partial differential
inequalities (PDIs for short) \cite{[BE],[BO1],[CL],[DZ2],[LPZ], dz,
[TY]}.

In the finite horizon framework Djehiche et al. have studied the
multi-modes switching problem in using probabilistic tools. For
general stochastic processes, they have shown that a value of the
problem and an optimal strategy exits. The partial differential
equation version of this work has been carried out by El-Asri and
Hamad\`ene \cite{[EH]}. They showed that when the price process
$X_t$ is solution of a Markovian standard differential equation,
then with this problem is associated a system of variational
inequalities with interconnected obstacles for which they provide a
solution in viscosity sense. This solution is bind to the value
function of the problem. The solution of the system is unique.

In the case when the horizon is infinite, there still much to do and
this is the novelty of this paper. Actually, authors treat mainly
the case when the price process $X_t$ is of Markovian It\^{o} type,
the switching costs are deterministic functions of time $t$ and the
profit functions are deterministic functions of $(t,X_t)$ and have
linear growth at most (see e.g. \cite{[BE], [BO1], [DZ2], [LPZ],
dz}). Therefore the main objective of this paper is to fill in the
gap between finite and infinite horizon by providing a complete
treatment of the optimal multiple switching problem in infinite
horizon when the price is only a continuous process. This is what we
did in the first part of this paper. Actually inspired by the work
of Djehiche et al. \cite{[DHP]}, using probabilistic tools such the
Snell envelope of processes and BSDEs we provide a verification
theorem which shapes the problem and then we have constructed a
solution for this latter. This solution provides an optimal strategy
for the switching problem. Later on, in the Markovian framework of
randomness, i.e. in the case when $X$ is a solution of a SDE, we
show that with the value function of the problem is associated an
uplet of deterministic functions $(v^1,\dots,v^m)$ which is the
unique solution of the following system of partial differential
inequalities (PDIs for short): \be\label{sysintro}\left\{
\begin{array}{l}
\min\{v_i(x)- \max\limits_{j\in{\cal
I}^{-i}}(-g_{ij}(x)+v_j(x)),rv_i(x)- {\cal A}v_i(x)-\psi_i(x)\}=0\\
\forall x\in \R^k,\,\,i\in{\cal I}=\{1,...,m\},
\end{array}\right.
\ee where $\cal A$ an infinitesimal generator associated with a
diffusion process and ${\cal I}^{-i}:=\{1,...,i-1,i+1,...,m\}$. This
system is the deterministic version of the Verification Theorem of
the optimal multi-modes switching problem in infinite horizon.


This paper is organized as follows: In Section 2, we formulate the
problem and we give the related definitions.  In Section 3, we
introduce the optimal switching problem under consideration and give
its probabilistic Verification Theorem. It is expressed by means of
a Snell envelope of processes. Then we introduce the approximating
scheme which enables to construct a solution for the Verification
Theorem. Moreover we give some properties of that solution,
especially the dynamic programming principle. Section 4 is devoted
to the connection between the optimal switching problem, the
Verification Theorem and the associated system of PDIs. This
connection is made through backward stochastic differential
equations with one reflecting obstacle in the case when randomness
comes from a solution of a standard stochastic differential
equation. Further we provide some estimate for the optimal strategy
of the switching problem which, in combination with the dynamic
programming principle, plays a crucial role in the proof of
existence of a solution for (\ref{sysintro}) which we address. In
Section 5, we show that the solution of PDIs is unique in the class
of continuous functions which satisfy a polynomial growth condition.
In Section 6, we give some numerical examples.$\Box$

\section{Assumptions and formulation of the problem}
Throughout this paper $k$ is a fixed integer positive constant. Let
us now consider the followings assumption:
\medskip

\indent $\bf H1$:  $b:R^k\rightarrow \R^{k}$ and
$\sigma:\R^k\rightarrow \R^{k\times d}$ are two continuous functions
for which there exists a constant $C\geq 0$ such that for any $x,
x'\in \R^k$\be \label{regbs1}|b(x)|+ |\sigma(x)|\leq C(1+|x|) \quad
\mbox{ and } \quad |\sigma(x)-\sigma(x')|+|b(x)-b(x')|\leq
C|x-x'|\ee

$\bf H2$:  for $i,j \in {\cal I}=\{1,...,m\}$,
$g_{ij}:\R^k\rightarrow \R$ is a continuous function. Moreover we
assume that there exists a constant $\alpha >0$ such that for any
$x\in
 \R^k$, \be  \frac{1}{\alpha}\leq g_{ij}(x)\leq \alpha,\quad \forall i,j \in {\cal I},\quad i
\neq j . \ee

$\bf H3$:   for $i\in {\cal I}$ $\psi_i:\R^k\rightarrow \R$ is a
continuous function of polynomial growth, $i.e.$, there exist a
constant $C$ and $\gamma$ such that for each $i\in \cal I$: \be
\label{polycond} |\psi_i(x)| \leq C(1+|x|^\gamma),\, \, \forall x\in
\R^k. \ee

\medskip

We now consider the following system of $m$ variational inequalities
with inter-connected obstacles:  $\forall \,\,i\in {\cal I}$ \be
\label{sysvi1}
\begin{array}{l}
\min\{v_i(x)- \max\limits_{j\in{\cal
I}^{-i}}(-g_{ij}(x)+v_j(x)),rv_i(x)- {\cal A}v_i(x)-\psi_i(x)\}=0
\end{array}.
\ee where ${\cal I}^{-i}:={\cal I}-\{i\}$, $r$ is a positive
discount factor and ${\cal A}$ is the following infinitesimal
generator:
\begin{equation}
\label{generateur} {\cal A}=\frac{1}{2}\sum_{i,j=1,k}(\sigma
\sigma^*)_{ij}(x)\frac{\partial^2} {\partial x_i \partial
x_j}+\sum_{i=1,k} b_i(x)\frac{\partial}{\partial x_i}\,;
\end{equation}
hereafter the superscript $(^*)$ stands for the transpose, $Tr$ is
the trace operator and finally $<x,y>$ is the inner product of
$x,y\in \R^k$.
\medskip

The main objective of this paper is to focus on the existence and
uniqueness of the solution in viscosity sense of (\ref{sysvi1})
whose definition is:

\begin{axiom} Let $(v_1,...,v_m)$ be a $m$-uplet of continuous functions defined on
$\R^k$, $\R$-valued. The $m$-uplet $(v_1,...,v_m)$ is called:
\begin{itemize}
\item [$(i)$] a viscosity supersolution (resp. subsolution) of the system (\ref{sysvi1})
if for each fixed $i\in {\cal I}$, for any $x_0\in \R^k$ and any
function $\varphi_i \in C^{1,2}(\R^k)$ such that
$\varphi_i(x_0)=v_i(x_0)$ and $x_0$ is a local maximum of $\varphi_i
-v_i$ (resp. minimum), we have: \be
\begin{array}{l}
\min\left\{v_i(x_0)- \max\limits_{j\in{\cal
I}^{-i}}(-g_{ij}(x_0)+v_j(x_0)),\right.\\\qquad\qquad \qquad
\left.r\varphi_i(x_0)-{\cal A}\varphi_i(x_0)-\psi_i(x_0)\right\}\geq
0 \quad (\mbox{resp.} \leq 0).
\end{array}
\ee
\item [$(ii)$] a viscosity solution if it is both a viscosity supersolution and
subsolution. $\Box$
\end{itemize}
\end{axiom}

There is an equivalent formulation of this definition (see e.g.
{\cite{[CIL]}) which we give since it will be useful later. So
firstly we define the notions of superjet and subjet of a continuous
function $v$.
\begin{axiom}
Let $v \in C(\R^k)$, $x$ an element of $ \R^k$ and finally $S_k$ the
set of $k \times k$ symmetric matrices. We denote by $J^{2,+} v(x)$
(resp. $J^{2,-} v(x)$), the superjets (resp. the subjets) of $v$ at
$x$, the set of pairs $(q,X)\in \R^k \times S_k$ such that:
$$\begin{array}{c}
v(y)\leq v(x) +\langle q,y-x \rangle +\frac{1}{2}\langle
X(y-x),y-x\rangle+o(|y-x|^2) \\
(resp.\quad v(y)\geq v(x) +\langle q,y-x \rangle +\frac{1}{2}\langle
X(y-x),y-x\rangle+o(|y-x|^2)). \Box
\end{array}$$
\end{axiom}
Note that if $\varphi-v$ has a local maximum (resp. minimum) at $x$,
then we obviously have:
$$\left(D_x \varphi(x),D^{2}_{xx}\varphi(x)\right)
\in J^{2,-} v (x) \,\,\, (\mbox{resp. } J^{2,+} v (x)). \Box$$

We now give an equivalent definition of a viscosity solution of the
elliptic system with inter-connected obstacles (\ref{sysvi1}).
\begin{axiom}Let $(v_1,...,v_m)$ be a $m$-uplet of continuous functions defined
on $\R^k$ and $\R$-valued. The $m$-uplet $(v_1,...,v_m)$ is called a
viscosity supersolution (resp. subsolution) of (\ref{sysvi1}) if for
any $i\in {\cal I}$, $x\in \R^k$ and $(q,X)\in J^{2,-} v_i (t,x)$
(resp. $J^{2,+} v_i (x)$),
$$ min \left\{v_i(x)-\max\limits_{j\in{\cal I}^{-i}}(-g_{ij}(x)+
v_j(x)),rv_i(x) -\frac{1}{2}Tr[\sigma^{*} X \sigma] -\langle b,q
\rangle-\psi_{i}(x)\right\}\geq 0 \,\,(resp. \leq 0).$$ It is called
a viscosity solution if it is both a viscosity subsolution and
supersolution .$\Box$
\end{axiom}

As pointed out previously we will show that system (\ref{sysvi1})
has a unique solution in viscosity sense. This system is the
deterministic version of the verification theorem of the optimal
$m$-states switching problem in infinite horizon which is well
documented in \cite{[DHP]} in the case of finite horizon and which
we will describe briefly in the next section.
\section{The  optimal $m$-states  switching problem}
\subsection{Setting of the problem}
Let $(\Omega, {\cal F}, P)$ be a fixed probability space on which is
defined a standard $d$-dimensional Brownian motion $B=(B_t)_{t\geq
0}$ whose natural filtration is $(\cF_t^0:=\sigma \{B_s, s\leq
t\})_{t\geq0}$. Let $ \bF=(\cF_t)_{ t\geq0}$ be the completed
filtration of $(\cF_t^0)_{t\geq0}$ with the $P$-null sets of ${\cal
F}$, hence $(\cF_t)_{t\geq0}$ satisfies the usual conditions,
$i.e.$, it is right continuous and complete. Furthermore, let:

- ${\cal P}$ be the $\sigma$-algebra on $[0,+\infty)\times \Omega$
of $\bF$-progressively measurable sets;

- ${\cal M}^{2,k}$ be the set of $\cal P$-measurable and
$\R^k$-valued processes $w=(w_t)_{t\geq0}$ such that
$E[\int_0^{+\infty}|w_s|^2ds]<\infty$  and ${\cal S}^2$  be the set
of $\cal P$-measurable, continuous processes ${w}=({w}_t)_{t\geq 0}$
such that $E[\sup_{t\geq0}|{w}_t|^2]<\infty$;

-  for any stopping time $\tau \in \R^+$, ${\cal T}_\tau$ denotes
the set of all stopping times $\theta$ such that $\tau \leq \theta;
$

-  for any stopping time $\tau$, $ {\cal F}_{\tau}$ is the
$\sigma$-algebra on $\Omega$ which contains the sets $A$ of
$\cal{F}$ such that $ A \cap \{\tau \leq t\}\in {\cal F}_t$ for
every $t\geq 0$. $\Box$
\medskip

A decision (strategy) of the problem of multiple switching, on the
one hand, consists of the choice of a sequence of nondecreasing
stopping times $(\tau_n)_{n\geq1}$ $(i.e. \tau_n \leq \tau_{n+1}$)
where the manager decides to switch the activity from its current
mode to another one. On the other hand, it consists of the choice of
the mode $\xi_n$, a r.v. ${\cal F}_{\tau_n}$-measurable with values
in ${\cal I}$, to which the production is switched at $\tau_n$.
Therefore the admissible management strategies are the pairs
$(\delta,\xi):=((\tau_n)_{n\geq 1},(\xi_n)_{n\geq 1})$ and we denote
by $\cal D$ the set of these strategies.


Let now $X:=(X_t)_{t\geq0}$ be an $\cal P$-measurable, $\R^k$-valued
continuous stochastic process which stands for the market price of
$k$ factors which determine the market price of the commodity. On
the other hand, assuming that the production activity is in mode 1
at the initial time $t = 0$, let $(u_t)_{t\geq0}$ denote the
indicator of the production activity's mode at time $t\in \R^+$ :
\begin{equation}
u_t=\ind_{[0,\tau_1]}(t)+\sum_{n\geq1}\xi_n
\ind_{(\tau_{n},\tau_{n+1}]}(t).
\end{equation}
Then for any $t\geq0$, the state of the whole economic system
related to the project at time $t$ is given by the vector:
\begin{equation}
\begin{array}{ll}
(t, X_t, u_t)\in \R^+\times \R^k \times {\cal I}.
\end{array}
\end{equation}

Finally, let $\psi_i(X_t)$ be the instantaneous profit when the
system is in state $(t, X_t, i)$, and for $i,j \in {\cal I} \quad
i\neq j$, let $g_{ij}(X_t)$ denote the switching cost of the
production at time $t$ from the current mode $i$ to another mode
$j$. When the plant is run under the strategy
$(\delta,\xi)=((\tau_n)_{n\geq 1},(\xi_n)_{n\geq 1})$ the expected
total profit is given by:
$$\begin{array}{l} J(\delta,\xi)=E[\integ{0}{+\infty}e^{-rs}\psi_{u_s}(X_s)ds
-\sum_{n\geq 1}
e^{-r\tau_n}g_{u_{\tau_{n-1}}u_{\tau_n}}(X_{\tau_{n}})].
\end{array}
$$
Then the problem we are interested in is to find an optimal
strategy, $i.e$, a strategy $(\delta^*,\xi^*)$ such that
$J(\delta^*,\xi^*)\ge J(\delta,\xi)$ for any $(\delta,\xi)\in \cal
D$.
\medskip

Note that in order that the quantity $J(\delta,\xi)$ makes sense we
assume throughout this paper that for any $i\in {\cal I}$ the
processes $(e^{-rt}\psi_i(X_t))_{t\geq0}$ belong to ${\cal
M}^{2,1}$. On the other hand there is a bijective correspondence
between the pairs $(\delta,\xi)$ and the pairs $(\delta,u)$. Then
throughout this paper one refers indifferently to $(\delta,\xi)$ or
$(\delta,u)$.
\subsection{ The Verification Theorem}

To tackle the problem described above in the finite horizon case,
Djehiche et al. \cite{[DHP]} have introduced a Verification Theorem
which is expressed by means of Snell envelope of processes which we
describe briefly now. The Snell envelope of a stochastic process
$(\eta_t)_{t\geq0}$ of ${\cal S}^2$ (with a possible positive jump
at $+\infty$ and $\lim\limits_{t\rightarrow\infty}\eta_t=M \in
L^2(\Omega, {\cal F}, P)$) is the lowest supermartingale
$R(\eta):=(R(\eta)_t)_{t\geq0}$ of ${\cal S}^2$ such that for any
$t\geq0$, $R(\eta)_t\geq \eta_t$. It has the following expression:
$$\forall t\geq0, R(\eta)_t=esssup_{\tau \in {\cal
T}_t}E[\eta_\tau|\bF_t] \quad \mbox{(then it satisfies
}\lim\limits_{t\rightarrow +\infty}R(\eta)_{t}=M.)$$ For more
details on the Snell envelope notion on can see e.g. \cite{[CK],
Elka, ham}.
\medskip

The Verification Theorem for the $m$-states optimal switching
problem in infinite horizon is the following:
\begin{theo}.\label{thmverif}
\noindent Assume that there exist $m$ processes
$(Y^i:=(Y^i_t)_{t\geq0}, i=1,...,m)$  of ${\cal S}^2$ such that: \be
\begin{array}{l}
\label{eqvt} \forall t\geq0,\,\,e^{-rt}Y^i_t=\esssup_{\tau \geq
t}E[\int_t^\tau e^{-rs}\psi_i(X_s)ds +e^{-r\tau}\max\limits_{j\in
{\cal I}^{-i}}(-g_{ij}(X_\tau)+Y^j_\tau)|\cF_t],\quad
\lim\limits_{t\rightarrow+\infty}(e^{-rt}Y^i_{t})=0.
\end{array}
\ee Then:
\begin{itemize}
\item[$(i)$]
$ Y^1_0=\sup \limits_{(\delta,\xi) \in {\cal D}}J(\delta,u). $
\item[$(ii)$] Define the sequence of $\bF$-stopping times $\delta^{*}=(\tau_n^*)_{n\geq 1}$ as follows :
$$
\begin{array}{lll}
\tau^*_{1}&=&\inf\{s\geq 0,\quad Y1_s=\max\limits_{j\in{{\cal I}^{-1}}}(-g_{1j}(X_s)+Y^j_s)\},\\
\tau^*_{n}&=&\inf\{s\geq \tau^*_{n-1},\quad
Y^{u_{\tau^*_{n-1}}}_s=\max\limits_{k\in {\cal I}\backslash
\{u_{\tau^*_{n-1}}\}} (-g_{u_{\tau^*_{n-1}}k}(X_s)+Y^k_s)\}, \quad
\mbox{for}\quad n\geq 2,
\end{array}
$$
where:
\begin{itemize}
\item[$\bullet$] $u_{\tau^*_1}=\sum \limits_{j\in {\cal I}} j \ind_ {\{\max\limits_{k\in {\cal I}^{-1}} (-g_{1k}(X_{\tau^*_1})+Y^k_{\tau^*_1})=-g_{1j}(X_{\tau^*_1})+Y^j_{\tau^*_1}\}};$
\item[$\bullet$] for any $n\geq1$ and $t\geq \tau^*_n,$
$Y^{u_{\tau^*_n}}_t =\sum\limits_{j\in {\cal I}}
\ind_{[u_{\tau^*_n}=j]}Y^j _t$
\item[$\bullet$] for any $n\geq 2, \,\,u_{\tau^*_n}=l$ on the set
$$\left\{\max\limits_{k\in {\cal I}
\backslash \{{ u_{\tau^*_{n-1}}}\}} (-g_{u_{\tau^*_{n-1}}
k}(X_{\tau^*_{n}})+ Y^k_{\tau^*_{n}})=-g_{u_{\tau^*_{n-1}
l}}(X_{\tau^*_n})+Y^l_{\tau^*_n}\right\}$$with \,
$g_{u_{\tau^*_{n-1} k}}(X_{\tau^*_n})= \sum\limits_{j\in {\cal
I}}\ind_{[u_{\tau^*_{n-1}}=j]}g_{j k}(X_{\tau^*_n})$ and ${\cal
I}\backslash \{u_{\tau^*_{n-1}}\}=\sum\limits_{j\in {\cal
I}}\ind_{[u_{\tau^*_{n-1}}=j]}{\cal I}^{-j}$.
\end{itemize}
Then the strategy $(\delta^*,u^*)$ satisfies $E[\sum_{n\geq0}
e^{-r\tau^*_{n}}]<+\infty$ and it is optimal i.e.
$J(\delta^*,u^*)\geq J(\delta,u)$ for any $(\delta,u)\in \cal D$.
 $\Box$
\end{itemize}
\end{theo}
\medskip
\no {\it Proof}. The arguments of this proof are standard, based on
the properties the Snell envelope. We defer the proof in the
Appendix.$\Box$

\medskip
The issue of existence of the processes $Y^1,...,Y^m$ which satisfy
(\ref{eqvt}) is also addressed in \cite{[DHP]}. For $n\geq 0$ let us
define the processes $(Y^{n,1},...,Y^{n,m})$ recursively as follows:
for $i\in {\cal I}$ we set,
\begin{equation}\label{y0}
e^{-rt}Y^{0,i}_t=E[\integ{t}{+\infty}e^{-rs}\psi_i(X_s)ds|{\cal
F}_t],\,\, t\geq0,
\end{equation}
and for $n\geq 1$,
\begin{equation}
\label{eq24} e^{-rt}Y^{n,i}_t=\mbox{ess sup}_{\tau\geq t}
E[\integ{t}{\tau}e^{-rs}\psi_i(X_s)ds+e^{-r\tau}\max\limits_{k\in
{\cal I}^{-i}}(-g_{ik}(X_\tau)+Y^{n-1,k}_\tau) |{\cal F}_t],\,\,
t\geq0.
\end{equation}
Then the sequence of processes $((Y^{n,1},...,Y^{n,m}))_{n\geq 0}$
have the following properties:
\begin{pro} (\cite{[DHP]}, Pro.3 and Th.2)
\begin{itemize}
\item[$(i)$] for any $i\in {\cal I}$ and
$n\geq 0$, the processes $Y^{n,1},...,Y^{n,m}$ are well-posed,
continuous and belong to ${\cal S}^2$, and verify
\be\label{croi}\forall t\geq0,\,\,e^{-rt}Y^{n,i}_t\leq
e^{-rt}Y^{n+1,i}_t\leq
E[\int_t^{+\infty}e^{-rs}\{\max_{i=1,m}|\psi_i(X_s)|\}ds|{\cal
F}_t];\ee

\item[$(ii)$] there exist $m$ processes $Y^1,...,Y^m$ of ${\cal S}^2$ such
that for any $i\in {\cal I}$:
\begin{itemize}
\item[$(a)$]  $\forall t\geq0$, $Y^i_t=\lim_{n\rightarrow
\infty}\nearrow Y^{n,i}_t$ 
\item[$(b)$] $\forall t\geq0$,
\begin{equation}
\label{eq26}e^{-rt}{Y}^{i}_t=\mbox{ess sup}_{\tau\geq
t}E[\integ{t}{\tau}e^{-rs}\psi_i(X_s)ds+ e^{-r\tau}\max\limits_{k\in
{\cal I}^{-i}}(-g_{ik}(X_\tau)+{Y}^{k}_\tau) |{\cal F}_t]
\end{equation}
i.e. ${Y}^1,...,{Y}^m$ satisfy the Verification Theorem
\ref{thmverif} ;
\item[$(c)$] $\forall t\geq0$,
\begin{equation}
\label{eq27}e^{-rt}{Y}^{i}_t=esssup_{(\delta,\xi)\in {\cal
D}^i_t}E[\integ{t}{+\infty}e^{-rs}\psi_{u_s}(X_s)ds -\sum_{n\geq 1}
e^{-r\tau_n}g_{u_{\tau_{n-1}}u_{\tau_n}}(X_{\tau_{n}}) |{\cal F}_t]
\end{equation}
where ${\cal D}^i_t=\{(\delta,\xi)=((\tau_n)_{n\geq
1},(\xi_n)_{n\geq 1}) \mbox { such that } u_0 =i \mbox{ and }
\tau_1\geq t \}$. This characterization means that if at time $t$
the production activity is in its regime $i$ then the optimal
expected profit is $Y^i_t$.
\item[$(d)$] the processes $Y^1,...,Y^m$ verify the dynamical programming principle of the
$m$-states optimal switching problem, $i.e.$, $\forall t\leq T$,
\begin{equation}\!\!\!\!\!
\label{dpp}
\begin{array}{ll}
e^{-rt}Y^i_t&=\esssup_{(\delta,u) \in {\cal D}_t^i}
E[\integ{t}{\tau_n}e^{-rs}\psi_{u_s}(X_s)ds -\sum_{1\leq k \leq n}
e^{-r\tau_k}g_{u_{\tau_{k -1}}u_{{\tau_k}}}(X_{{\tau}_{k}})
+e^{-r\tau_n}Y^{u_{\tau_n}}_{\tau_n}|{\cal F}_t].\Box
\end{array}
\end{equation}
\end{itemize}
\end{itemize}
\end{pro}

Note that except $(ii-d)$, the proofs of the other points are the
same as in \cite{[DHP]} in the framework of finite horizon. The
proof of $(ii-d)$ can be easily deduced in using relation
(\ref{eq26}). Actually from (\ref{eq26}) for any $i\in {\cal I}$,
$t\geq 0$ and $(\delta,\xi)\in {\cal D}^i_t$ we have: \be
\label{eq28}e^{-rt}Y^i_t\geq
E[\integ{t}{\tau_n}e^{-rs}\psi_{u_s}(X_s)ds -\sum_{1\leq k \leq n}
e^{-r\tau_k}g_{u_{\tau_{k -1}}u_{u_{\tau_k}}}(X_{{\tau}_{k}})
+e^{-r\tau_n
 }Y^{u_{\tau_n}}_{\tau_n}|{\cal
F}_t]. \ee Next using the optimal strategy we obtain the equality
instead of inequality in (\ref{eq28}). Therefore the relation
(\ref{dpp}) holds true. $\Box$
\begin{rem} \label{unic}The characterization (\ref{eq27}) implies that the processes $Y^1,...,Y^m$ of ${\cal S}^2$ which satisfy the Verification
Theorem are unique.
\end{rem}
\section{Existence of a solution for the system of variational inequalities}
\subsection{Connection with BSDEs with  one reflecting barrier}
Let $x\in \R^k$ and let $X^{x}$ be the solution of the following
standard SDE:
\begin{equation}\label{sde}
dX_t^x=b(X^x_t)dt+\sigma(X^x_t)dB_t, \quad X^x_0=x
\end{equation}where the functions $b$ and $\sigma$ are the ones of
$\bf H1$. These properties of $\sigma$ and $b$ imply in particular
that $X^{x}$ solution of the standard SDE (\ref{sde}) exists and is
unique in $\R^k $. The operator $\cal A$  defined in
(\ref{generateur}) is the infinitesimal generator associated with
$X^{x}$.

 In the following result we collect some properties of
$X^{x}$.

\bp \label{estimx}(see e.g. \cite{[RY]})  The process $X^{x}$
satisfies the following estimates:
\begin{itemize}
\item [$(i)$] For any $q\geq 2$ there exists $C_q$ such that,
\begin{equation}\label{estimat1}
E[|X^{x}_t|^q]\leq C_qe^{C_qt}(1+|x|^q)\quad \forall t\geq 0.
\end{equation}
\item[$(ii)$] There exists a constant $C$ such that for any  $x,x'\in \R^k$ and $T\geq 0$,
\begin{equation}\label{estimat3}
E[\sup\limits_{0\leq s \leq  T}|X^{x}_s-X^{x'}_s|^2]\leq
Ce^{CT}|x-x'|^2. \Box
\end{equation}
\end{itemize}
\ep

In the sequel we consider the following condition:
\medskip

\noindent $\bf H4$: Assume $ \gamma\geq 2$ and
\begin{equation}\label{discount}-r+C_{\gamma}<0,\end{equation} where $\gamma$ is the growth exponent of
the functions $\psi_i$ and $C_{\gamma}$ is the constant in
(\ref{estimat1}). $\Box$
\begin{rem}: \label{unic}If  $\gamma<2$, there exists a constant $\gamma_1\geq2$ such that $\gamma_1$
 verifies  the growth exponent of the functions $\psi_i$.
\end{rem}

We are going now to introduce the notion of a BSDE with one
reflecting barrier considered in \cite{[HLW]}. This notion will
allow us to make the connection between the variational inequalities
system (\ref{sysvi1}) and the $m$-states optimal switching problem
described in the previous section.
\medskip

Let us introduce the pair of  process $(Y^{x},Z^{x})\in {\cal S}^2
\times {\cal M}^{2,d}$ solution of the following BSDE:
\begin{equation}\label{bsde,sde}
Y^{x}_s=Y_T^x
+\int_{s}^{T}F(X_l^{x},Y_l^{x},Z_l^{x})dl-\int_{s}^{T}Z^{x}_ldB_l,\quad
\mbox{for all}\quad T\geq0 \quad \mbox{and}\quad t\leq T,
\end{equation}
where $F :\R^k \times \R \times \R^d\rightarrow \R$ is continuous
and satisfies: there exist a continuous increasing function $\phi :
\R^+ \rightarrow \R^+$ and constant $K$, $K'$, $\mu<0$, $p>0$ such
that,
\begin{equation}
\label{assumption-pardoux}
\begin{array}{lll}
|F (x,y,z)| \leq K'(1+|x|^p +\phi(|y|)+|z|),\\
\langle y-y',F(x,y,z)-F(x,y',z)\rangle \leq \mu|y-y'|^2,\\
|F (x,y,z)-F(x,y,z')|\leq K||z-z'||.
\end{array}
\end{equation}
We assume moreover that for some $\lambda>2\mu + K2,$
\begin{equation}
\label{assumption-pardoux1} E[\int_{0}^{+\infty}e^{\lambda s}|F
(X_s^{x},0,0)|^2 ds]<+\infty,
\end{equation}
which essentially implies that $\lambda +C_{2\gamma} <0$.\\

Let us consider the following semilinear elliptic PDE in $\R^k$:
\begin{equation}\label{viscosity}
{\cal A}u(x)+F(x,u(x),\sigma (x)^*\nabla u(x))=0,\quad x\in \R^k.
\end{equation}

 Then we have the following result:
\begin{theo}(\cite{[PA]}, Th. 5.2)\label{pardoux} Under the above assumptions,
$u(x)=Y^x_0$ is a continuous function and it is a viscosity solution
of (\ref{viscosity}) which satisfies,
\begin{equation}
\label{assumption-pardoux2} |Y^{x}_0|\leq C
\sqrt{E[\int_{0}^{+\infty}e^{\lambda s}|F (X_s^{x},0,0)|^2 ds}],
\end{equation}
for any $\lambda>2\mu + K2.$$\Box$

\end{theo}
 Let us now introduce the following functions:
\begin{itemize}
\item [$(i)$]
$f:\R^k\rightarrow \R$ is continuous and of polynomial growth,
$i.e.$, there exist some positive constants $C$ and $\gamma$ such
that: \be \label{polycond} |f(x)| \leq C(1+|x|^{\gamma}),\, \,
\forall x\in \R^k. \ee
\item[$(ii)$] $h:\R^{k}\rightarrow \R$ is
 continuous and bounded.
 \end{itemize}
  Then we have the following
result related to BSDEs with one reflecting barrier:
\begin{theo} For any $x\in \R^k$, there exits a unique triple
of processes $(Y^{x},Z^{x},K^{x})$ such that:
\begin{equation}\label{BSDE}\left\{
\begin{array}{l}
Y^{x}, K^{x}\in {\cal S}^2 \mbox{ and }Z^{x}\in {\cal
M}^{2,d};\,K^{x}
\mbox{ is  non-decreasing and }K^{x}_0=0,\\
e^{-rs}Y^{x}_s=\int_{s}^{+\infty}e^{-rl}f(X_l^{x})dl-\int_{s}^{+\infty}Z^{x}_ldB_l+
K_{+\infty}^{x}-K^{x}_s, \,\, \\
e^{-rs}Y^{x}_s\geq e^{-rs}h(X^{x}_s),\, \forall s \geq 0\mbox{ and }
\int_{0}^{+\infty}(e^{-rl}Y^{x}_l-e^{-rl}h(X^{x}_l))dK^{x}_l=0.
\end{array}
\right. \end{equation} Moreover the following characterization of
$Y^{x}$ as a Snell envelope holds true: \be \label{snellenv}\forall
s\geq 0,\,\,e^{-rs}Y^{x}_s=esssup_{\tau \in {\cal T}_s}E[\int_s^\tau
e^{-rl}f(X_l^{x})dl +e^{-r\tau}h(X^{x}_\tau)|{\cal F}_s].\ee

On the other hand there exists a deterministic continuous with
polynomial growth function $u:\R^{k}\rightarrow \R$ such
that:$$\forall x\in\R^k \quad Y^{x}_0=u(x).$$ Moreover the function
$u$ is the viscosity solution in the class of continuous function
with polynomial growth of the following PDE with
obstacle:\begin{equation}\label{viscosity-equation}
\begin{array}{l}
\min\{u(x)- h(x), ru(x)-{\cal A}u(x)-f(x)\}=0.
\end{array}
\end{equation}
\end{theo}
$Proof$: Existence and uniqueness of the triple
$(Y^x_t,Z^x_t,K^x_t)_{t\geq0}$ of (\ref{BSDE}) follow from Theorem
3.2 in \cite{[HLW]}. Now we consider the infinite horizon BSDE:
\begin{equation}\label{penalization}
 ^nY^{x}_se^{-rs}=\int_{s}^{+\infty}e^{-rl}f(X_l^{x})dl-\int_{s}^{+\infty}Z^{n,x}_l dB_l+ \int_{s}^{+\infty}ne^{-rl}(^n Y^{x}_l
 -h(X^{x}_l))^{-}dl.
\end{equation}
From Theorem 1 in \cite{[CH]} there exists a unique solution
$(^nY^{x},Z^{n,x})\in{\cal S}^2\times {\cal M}^{2,d}$ satisfying the
BSDE
(\ref{penalization}).\\
\indent Next let us define
$$K^{n,x}_s=\int_{0}^{s}ne^{-rl}(^nY^{x}_l-h(X^{x}_l))^-dl,$$
then
$$
\begin{array}{ll}
\int_{0}^{+\infty}e^{-rl}(^nY^{x}_l-h(X^{x}_l)\wedge^nY^{x}_l)dK_l^{n,x}&=n\int_{0}^{+\infty}e^{-rl}(^nY^{x}_l-h(X^{x}_l)\wedge^nY^{x}_l)
e^{-rl}(^nY^{x}_l-h(X^{x}_l))^-dl\\&=0.
\end{array}
$$
Since $K^{n,x}$  is  non-decreasing and $K^{n,x}_0=0$, we rewrite
Eq. (\ref{penalization}) in RBSDE form
\begin{equation}\label{bsde-pena}\left\{
\begin{array}{l}
^nY^{x}_se^{-rs}=\int_{s}^{+\infty}e^{-rl}f(X_l^{x})dl-\int_{s}^{+\infty}Z^{n,x}_ldB_l+
K_{\infty}^{n,x}-K^{n,x}_s, \,\, \\
^nY^{x}_se^{-rs}\geq e^{-rs}(h(X^{x}_s)\wedge ^nY^{x}_s) ,\, \forall
s \geq 0\mbox{ and }
\int_{0}^{+\infty}e^{-rl}(^nY^{x}_l-h(X^{x}_l)\wedge
^nY^{x}_l)dK^{x}_l=0.
\end{array}
\right. \end{equation}  Then from property (\ref{snellenv}) we have:
\begin{equation}\label{envelope-penalisation}
^nY^{x}_se^{-rs}=esssup_{\tau \in {\cal T}_s}E[\int_s^\tau
e^{-rl}f(X_l^{x})dl +e^{-r\tau}(^nY^{x}_{\tau}\wedge
h(X^{x}_\tau))|{\cal F}_s].
\end{equation}
Note that if we define
$$f_n(t,x,y,z)=e^{-rt}f(x,y,z)+ne^{-rt}(y-h(x))^-$$
$$f_n(t,x,y,z)\leq f_{n+1}(t,x,y,z).$$
Then it follows from the comparison Theorem 2.2 in \cite{[HLW]}
$^nY_s^{x}e^{-rs}\leq ^{n+1}Y_s^{x}e^{-rs},$ $s\geq 0,$ a.s. and
from (\ref{snellenv}) and (\ref{envelope-penalisation})\quad$
^nY_s^{x}e^{-rs}\leq Y_s^{x}e^{-rs}.$  This implies that there exits
a c\`adl\`ag process $(\widetilde{Y}^x_s)_{s\geq 0}$ such that
$P-a.s$. for any $s\geq 0$,
$$^nY_s^{x}e^{-rs}\uparrow e^{-rs}\widetilde{Y}_s^{x}, \quad\quad a.s.$$
\noindent Let us actually show that $\tilde{Y}^x$ is \cadlag. By
(\ref{envelope-penalisation}), for any $n\geq 1$, the process
$(^nY^{x}_t+\int_{0}^{t}e^{-rs}f(X_s^x)ds)_{t\geq 0}$ is an
$\bF$-supermartingale which converges increasingly and pointwisely
to $(\widetilde{Y}^{x}_t+\int_{0}^{t}e^{-rs}f(X_s^x)ds)_{t\geq 0}$.
Therefore, the limit is also a \cadlag $\bF$-supermartingale (see
e.g. Dellacherie and Meyer (1980), pp. 86). Hence, the process
$\tilde{Y}^x$ is \cadlag.

Then it follows from Proposition 2 in \cite{[DHP]}, as $n\rightarrow
+\infty$,
\begin{equation}\label{envelope}
\widetilde{Y}^{x}_se^{-rs}=esssup_{\tau \in {\cal T}_s}E[\int_s^\tau
e^{-rl}f(X_l^{x})dl +e^{-r\tau}(\widetilde{Y}^{x}_{\tau}\wedge
h(X^{x}_\tau))|{\cal F}_s].
\end{equation}
 From(\ref{penalization}) we have:
$$
\begin{array}{ll}
E[\int_{s}^{+\infty}e^{-rl}(^n Y^{x}_l -h(X^{x}_l))^{-}dl]&= \frac{1}{n}E[^nY^{x}_se^{-rs}+\int_{s}^{+\infty}e^{-rl}f(X_l^{x})dl]\\
&\leq
\frac{1}{n}E[|Y^{x}_se^{-rs}|+\int_{s}^{+\infty}|e^{-rl}f(X_l^{x})|dl]\\&\leq
\frac{1}{n}(E[|Y^{x}_se^{-rs}|]+C\int_{s}^{+\infty}e^{-rl}e^{C_{\gamma}l}|x|^{\gamma}dl)
\end{array}
$$
for a constant $C$ independent of $n$ and 
$\bf H4$. Then
$$ E [\int_{s}^{+\infty}e^{-rl}(^n Y^{x}_l-h(X^{x}_l))^{-}dl] \leq
\frac{C_x}{n}.$$ Hence as $n\rightarrow +\infty$ we obtain,
$E[\int_{s}^{+\infty}e^{-rl}(\widetilde{Y}^{x}_l
-h(X^{x}_l))^{-}dl]=0$, and since $(\widetilde{Y}^x_s)_{s\geq 0}$
(resp. $h(x)$) is a c\`adl\`ag process (resp. continuous), we have
\begin{equation}\label{comparaison}
\widetilde{Y}^{x}_t \geq h(X^{x}_t). \end{equation} From
(\ref{snellenv}), (\ref{envelope}) and (\ref{comparaison}) we get:

$$\widetilde{Y}^x_t=Y^x_t \quad \forall t\geq 0.$$ Now rewrite Eq.
(\ref{penalization}) in differential form
$$
\begin{array}{l}
d(^nY^{x}_se^{-rs})=-[e^{-rs}f(X_s^{x})+ne^{-rs}(^nY^{x}_s -h(X^{x}_s))^{-}]ds+Z^{n,x}_s dB_s.\\
\end{array} 
$$
So for arbitrary $T> 0$ and $0\leq s\leq T$, Eq.
(\ref{penalization}) is equivalent to
\begin{equation}
\label{mode-penalization}
\begin{array}{l}
^nY^{x}_s=^nY^{x}_T+\int_{s}^{T}[(f(X_l^{x})+n(^nY^{x}_l -h(X^{x}_l))^{-})-r^nY^{x}_l]dl-\int_{s}^{T}\widetilde{Z}^{n,x}_l dB_l,\\
\end{array} 
\end{equation}
with $\widetilde{Z}^{n,x}_s= Z^{n,x}_se^{rs}$. Let us set $F_n(x,y,z)=f(x)+n(y -h(x))^{-})-ry$.\\
In order that it satisfies the assumptions of Theorem \ref{pardoux},
we just need to verify that $F_n$ satisfy condition
(\ref{assumption-pardoux}) and (\ref{assumption-pardoux1}). It is
obvious that $F_n$ satisfy (\ref{assumption-pardoux})
where $\mu >-r$, and we show that $F_n$ satisfy (\ref{assumption-pardoux1}).\\
From the polynomial growth of $f$ and since $h$ bounded and estimate
(\ref{estimat1}), we deduce
$$
\begin{array}{lll}
E[\int_{0}^{+\infty}e^{\lambda s}|F_n (X_s^{x},0,0)|^2 ds]&=E[\int_{0}^{+\infty}e^{\lambda s}|f(X_s^{x})+n(-h(X_s^{x}))^-|^2 ds]\\
&\leq 2 E[\int_{0}^{+\infty}e^{\lambda s}((1+|X_s^{x}|^\gamma)2+n^2C2) ds]\\
&\leq C\int_{0}^{+\infty}e^{\lambda s}e^{C_{2\gamma}
s}(|x|^{2\gamma}+n2) ds,
\end{array}
$$
for $\lambda+ C_{2\gamma}<0$. This proves assumption
(\ref{assumption-pardoux1}). Then
$$u_n(x)=^nY^{x}_0,$$and is a viscosity solution of the elliptic PDE
$${\cal A}u_n(x)+F_n(x,u_n(x),\sigma (x)^*\nabla u_n(x))=0.$$
\indent We now define $$u(x)=Y^{x}_0,\quad \forall x \in \R^{k},$$
which is a deterministic quantity. Let us admit for a moment the
following Lemma:
\begin{lem}\label{lemma}
The function $u$ is continuous in $R^k$.$\Box$
\end{lem}

%
From the previous results we have, for each $x\in \R^k ,$
$$u_n(x)\uparrow u(x)\quad \mbox{as}\quad n\rightarrow +\infty.$$
Since $u_n$ and $u$ are continuous, it follow from Dini's theorem that the above convergence is uniform on compacts.\\
 \indent We now show that $u$ is a subsolution of (\ref{viscosity-equation}). Let $x$ be a point at which $u(x)>h(x),$ and let $(q,X)\in
 J^{2,+}u(x).$
From Lemma 6.1 in \cite{[CIL]}, there exists sequences:
$$
\begin{array}{l}
n_j \rightarrow+\infty,\quad x_j\rightarrow x, \quad (q_j,X_j)\in
J^{2,+}u_{n_j}(x_j),
\end{array}
$$
such that $$(q_j,X_j)\rightarrow (q,X).$$ But for any $j$,
$$
\begin{array}{ll}
-\frac{1}{2}Tr[\sigma^{*} X_j \sigma] -\langle b,q_j \rangle - F_n(x_j,u_{n_j}(x_j),\sigma (x_j)^*\nabla u_{n_j}(x_j))\leq 0,\\
-\frac{1}{2}Tr[\sigma^{*} X_j \sigma] -\langle b,q_j \rangle -
f(x_j)-n_j(u_{n_j}(x_j) -h(x_j))^{-})+ru_{n_j}(x_j)\leq 0.
\end{array}
$$
From the assumption that $u(x)>h(x)$ and the uniform convergence of
$u_n,$ it follows that for $j$ large enough $u_{n_j}(x_j)>h(x_j)$.
Hence, taking the limit as $j\rightarrow +\infty$ in the above
inequality yields:
$$-\frac{1}{2}Tr[\sigma^{*} X \sigma] -\langle b,q \rangle -f(x)+ ru(x) \leq 0, $$
and we have proved that $u$ is a subsolution of (\ref{viscosity-equation}).\\\\
We now show that $u$ is a supersolution of
(\ref{viscosity-equation}). Let $x$ be arbitrary in $\R^k$, and
$(q,X)\in J^{2,-}u(x).$ We already know that $u(x)\geq h(x).$ By the
same argument as above, there exist sequences:
$$
\begin{array}{l}
n_j \rightarrow+\infty,\quad x_j\rightarrow x,\quad (q_j,X_j)\in
J^{2,-}u_{n_j}(x_j),
\end{array}
$$
such that $$(q_j,X_j)\rightarrow (q,X).$$ But for any $j$,
$$
\begin{array}{ll}
-\frac{1}{2}Tr[\sigma^{*} X_j \sigma] -\langle b,q_j \rangle - F_n(x_j,u_{n_j}(x_j),\sigma (x_j,i)^*\nabla u_{n_j}(x_j)))\geq 0,\\
-\frac{1}{2}Tr[\sigma^{*} X_j \sigma] -\langle b,q_j \rangle -
f(x_j)-n_j(u_{n_j}(x_j) -h(x_j))^{-})+ru_{n_j}(x_j)\geq 0.
\end{array}
$$
Hence,
$$-\frac{1}{2}Tr[\sigma^{*} X_j \sigma] -\langle b,q_j \rangle -f(x_j)+ ru_{n_j}(x_j) \geq 0, $$
and taking the limit as $j\rightarrow +\infty$, we conclude that:
$$-\frac{1}{2}Tr[\sigma^{*} X \sigma] -\langle b,q \rangle -f(x)+ ru(x) \geq 0. $$
 We conclude by showing that $u$ is of polynomial growth.
From (\ref{snellenv}) we have,

\begin{equation}\label{growth_pol}
\begin{array}{ll}
|Y^{x}_0|&\leq  sup_{\tau \geq 0}E[\int_0^\tau e^{-rs}|f(X_s^{x})|ds
+ |h(X^{x}_\tau)|\ind_{[\tau<+\infty]}]\\
&\leq sup_{\tau \geq 0}E[\int_0^\tau e^{-rs}|f(X_s^{x})|ds
+e^{-r\tau}|h(X^{x}_\tau)|]\\
&\leq E[\int_0^{+\infty} e^{-rs}|f(X_s^{x})|ds] +C_1.
\end{array}
\end{equation}
From polynomial growth of $f$ and $u(x)=Y^{x}_0$, we deduce that $u$
is of polynomial growth.
Now we proceed to the proof of Lemme1.\\
\indent $Proof$ of Lemma 2. It suffices to show that whenever
$x_n\rightarrow x$, $|Y_0^{x_n}-Y_0^x|\rightarrow0$.\\
From (\ref{snellenv}) we have,
$$Y^{x}_0=\sup_{\tau \in {\cal T}_0}E[\int_0^\tau
e^{-rl}f(X_l^{x})dl +e^{-r\tau}h(X^{x}_\tau)],$$
$$Y^{x_n}_0=\sup_{\tau \in {\cal T}_0}E[\int_0^\tau
e^{-rl}f(X_l^{x_n})dl +e^{-r\tau}h(X^{x_n}_\tau)]$$ then,
\begin{equation}\label{preuve_cont}
\begin{array}{ll}
|Y^{x_n}_0 -Y^{x}_0|&\leq \sup\limits_{\tau \in {\cal
T}_0}E[\int_0^\tau e^{-rl}|f(X_l^{x_n})-f(X_l^{x})|dl
+e^{-r\tau}|h(X^{x_n}_\tau)-h(X^{x}_\tau)|]\\
&\leq E[\int_0^{+\infty} e^{-rl}|f(X_l^{x_n})-f(X_l^{x})|dl]
+E[\sup\limits_{t\geq0}e^{-rt}|h(X^{x_n}_t)-h(X^{x}_t)|].
\end{array}
\end{equation}
In the right-hand side of (\ref{preuve_cont}) the first term
converges to 0 as $x_n\rightarrow x$. Next let us show that,
$$E[\sup\limits_{t\geq0}e^{-rt}|h(X^{x_n}_t)-h(X^{x}_t)|]\rightarrow0 \quad \mbox{as}\quad x_n\rightarrow x.$$
For any $T\geq 0$ we have
$$
E[\sup\limits_{t\geq0}e^{-rt}|h(X^{x_n}_t)-h(X^{x}_t)|]\leq
E[\sup\limits_{0\leq t\leq
T}e^{-rt}|h(X^{x_n}_t)-h(X^{x}_t)|]+E[\sup\limits_{t\geq
T}e^{-rt}|h(X^{x_n}_t)-h(X^{x}_t)|].$$ Since $h$ is bounded there
exists $C$ such that,
$$
E[\sup\limits_{t\geq0}e^{-rt}|h(X^{x_n}_t)-h(X^{x}_t)|] \leq
E[\sup\limits_{0\leq t\leq T}e^{-rt}|h(X^{x_n}_t)-h(X^{x}_t)|]+
Ce^{-rT}.
$$
For any $\rho> 0$ we have:
$$
\begin{array}{ll}\label{separa}
E[\sup\limits_{0\leq t\leq
T}e^{-rt}|h(X^{x_n}_t)-h(X^{x}_t)|]&=E[\sup\limits_{0\leq t\leq
T}e^{-rt}|h(X^{x_n}_t)-h(X^{x}_t)|\ind_{[\sup\limits_{t\leq T}
|X_t^{x_n}|+ \sup\limits_{t\leq T}|X_t^{x}|
\leq\rho]}]\\
&+E[\sup\limits_{0\leq t\leq
T}e^{-rt}|h(X^{x_n}_t)-h(X^{x}_t)|\ind_{[\sup\limits_{t\leq
T}|X_t^{x_n}|+\sup\limits_{t\leq T}|X_t^{x}|>\rho]}].
\end{array}
$$
But since $h$ is continuous then it is uniformly continuous on
compact subsets, then there exists $\pi:R^k\rightarrow R$ increasing
with $\pi(0)=0$, such that:
$$|h(X^{x_n}_t)-h(X^{x}_t)|\leq \pi(|X^{x_n}_t-X^{x}_t|),$$
we have
$$
\begin{array}{ll}
E[\sup\limits_{0\leq t\leq
T}e^{-rt}|h(X^{x_n}_t)-h(X^{x}_t)|\ind_{[\sup\limits_{t\leq T}
|X_t^{x_n}|+ \sup\limits_{t\leq T}|X_t^{x}| \leq\rho]}]&\leq
E[\sup\limits_{0\leq t\leq
T}\pi(|X^{x_n}_t-X^{x}_t|)\ind_{[\sup\limits_{t\leq T} |X_t^{x_n}|+
\sup\limits_{t\leq T}|X_t^{x}| \leq\rho]}]\\&\leq
E[\pi(\sup\limits_{0\leq t\leq
T}|X^{x_n}_t-X^{x}_t|)\ind_{[\sup\limits_{t\leq T} |X_t^{x_n}|+
\sup\limits_{t\leq T}|X_t^{x}| \leq\rho]}].
\end{array}$$
 Using the
continuity proprety (\ref{estimat3}), $\pi(0)=0$ and the Lebesgue
dominated convergence theorem to obtain that
\begin{equation} E[\sup\limits_{0\leq t\leq
T}e^{-rt}|h(X^{x_n}_t)-h(X^{x}_t)|\ind_{[\sup\limits_{t\leq T}
|X_t^{x_n}|+ \sup\limits_{t\leq T}|X_t^{x}|
\leq\rho]}]\rightarrow0\quad \mbox{as}\quad x_n\rightarrow x.
\end{equation}
 The second term satisfies:
$$
\begin{array}{ll}
E[\sup\limits_{0\leq t\leq
T}e^{-rt}|h(X^{x_n}_t)-h(X^{x}_t)|\ind_{[\sup\limits_{t\leq T}
|X_t^{x_n}|+ \sup\limits_{t\leq T}|X_t^{x}|>\rho]}]
\\
{}\qquad  \leq E[\sup\limits_{0\leq t\leq
T}e^{-2rt}|h(X^{x_n}_t)-h(X^{x}_t)|^2]\}^{\frac{1}{2}}\{E[\ind_{[\sup\limits_{t\leq
T} |X_t^{x_n}|+ \sup\limits_{t\leq T}|X_t^{x}|
> >\rho]}]\}^{\frac{1}{2}} \\ {}\qquad \leq E[\{\sup\limits_{0\leq t\leq
T}e^{-2rt}|h(X^{x_n}_t)-h(X^{x}_t)|^2]\}^{\frac{1}{2}}\{\
\rho^{-1}E[\sup\limits_{t\leq T} |X_t^{x_n}|+ \sup\limits_{t\leq
T}|X_t^{x}|]\}^{\frac{1}{2}}.
\end{array}
$$
Since $h$ is bounded, it follows that, when $x_n\rightarrow x$, the
right-hand side of the last inequality is smaller than
$\rho^{-\frac{1}{2}}C_x$. However, from previous results we have,
$$\limsup\limits_{x_n\rightarrow x}E[\sup\limits_{t\geq0}e^{-rt}|h(X^{x_n}_t)-h(X^{x}_t)|]\leq \rho^{-\frac{1}{2}}C_x + Ce^{-rT}.$$

As $\rho$ and $T$ are arbitrary then making $\rho\rightarrow
+\infty$ and $T\rightarrow +\infty$ to obtain that,
\begin{equation}\label{ineg1}
\lim\limits_{x_n\rightarrow
x}E[\sup\limits_{t\geq0}e^{-rt}|h(X^{x_n}_t)-h(X^{x}_t)|]=0.
\end{equation}
From (\ref{preuve_cont}) and (\ref{ineg1}), we deduce
$$|Y_0^{x_n}-Y_0^x|\rightarrow0 \quad\mbox{as}\quad
x_n\rightarrow x.\Box$$

\subsection{Existence of a solution for the system of variational inequalities}

Let $(Y^{1,x}_s,...,Y^{m,x}_s)_{s\geq 0}$ be the processes which
satisfy the Verification Theorem \ref{thmverif} in the case when the
process $X\equiv X^{x}$. Therefore using the characterization
(\ref{snellenv}), there exist processes $K^{i,x}$ and $Z^{i,x}$,
such that the triples ($Y^{i,x}, Z^{i,x},K^{i,x})$ are unique
solutions (thanks to Remark \ref{unic}) of the following reflected
BSDEs: for any $i=1,...,m$ we have,
\begin{equation}\left\{
\begin{array}{l}
Y^{i,x}, K^{i,x}\in {\cal S}^2 \mbox{ and }Z^{i,x}\in {\cal
M}^{2,d};\,K^{i,x}
\mbox{ is  non-decreasing and }K^{i,x}_0=0,\\
e^{-rs}Y^{i,x}_s=\int_{s}^{+\infty}e^{-rl}\psi_{i}(X_l^{x})ds-\int_{s}^{+\infty}Z^{i,x}_ldB_l+K_{+\infty}^{i,x}-K^{i,x}_s,
\,\,\,
s\in \R^+ ,\,\,\lim\limits_{s\rightarrow +\infty}(e^{-rs}Y^{i,x}_s)=0,\\
e^{-rs}Y^{i,x}_s\geq -e^{-rs}\max\limits_{j\in{\cal I}^{-i}}(-g_{ij}(X_s^{x})+Y^{j,x}_s),\,\, s\in \R^+,\\
\int_{0}^{+\infty}e^{-rl}(Y^{i,x}_l-\max\limits_{j\in{\cal
I}^{-i}}(-g_{ij}(X_l^{x})+Y^{j,x}_l))dK^{i,x}_l=0.
\end{array}
\right. \end{equation} Moreover we have the following result.
\begin{pro}There are deterministic functions $v^1,...,v^m$
$:\R^k\rightarrow \R$ such that:
$$\forall x\in \R^k,
Y_0^{i,x}=v^i(x), \,\,i=1,...,m.$$ Moreover the functions $v^i$,
$i=1,...,m,$ are of polynomial growth.
\end{pro}
$Proof$: For $n\geq 0$ let $(Y^{n,1,x}_s,...,Y^{n,m,x}_s)_{s\geq0}$
be the processes constructed in (\ref{y0})-(\ref{eq24}). Therefore
using an induction argument and Theorem 2 there exist deterministic
continuous with polynomial growth functions $v^{n,i}$ ($i=1,...,m$)
such that for any $x\in \R^k$,
 $Y^{n,i,x}_0=v^{n,i}(x)$. Using now
inequality (\ref{croi}) we get:
$$Y^{n,i,x}_t\le Y^{n+1,i,x}_t\leq
CE[\int_0^{+\infty}\{\max_{i=1,m}|e^{-rs}\psi_i(X^{x}_s)|\}ds]$$
since $Y^{n,i,x}_t$ is deterministic. Therefore combining the
polynomial growth of $\psi_i$ and estimate (\ref{estimat1}) for
$X^{x}$ we obtain:
$$v^{n,i}(x)\leq v^{n+1,i}(x)\leq C(1+|x|^\gamma)$$for a constant
$C$ independent of $n$. In order to complete the proof it is enough
now to set $v^i(x):=\lim_{n\rightarrow \infty}v^{n,i}(x), x\in \R^k$
since $Y^{n,i,x}\nearrow Y^{i,x}$ as $n \rightarrow \infty$. $\Box$
\medskip

We are now going to focus on the continuity of the functions
$v^1,...,v^m$. But first let us deal with some properties of the
optimal strategy which exist thanks to Theorem 1.

\bp \label{optimal-s} Let $(\delta,u)=((\tau_n)_{n\geq
1},(\xi_n)_{n\geq 1})$ be an optimal strategy, then there exists a
constant $C$ which does not depend on $t$ and $x$ such that: \be
\label{estiopt}\forall n\geq 1,\,\,E[e^{-r\tau_n}]\leq \frac
{C(1+|x|^\gamma)}{n}.\ee \ep $Proof$: Recall the characterization of
(\ref{eq27}) that reads as:
$$
\begin{array}{l}
Y^{i,x}_0=sup_{(\delta,u) \in {\cal
D}}E[\int_0^{+\infty}e^{-rs}\psi_{u_s}(X_s^{x})ds-\sum_{k\geq 1}
e^{-r\tau_k}g_{u_{\tau_{k-1}}u_{\tau_k}}(X^{x}_{\tau_{k}})].
\end{array}
$$
Now if $(\delta,u)=((\tau_n)_{n\geq 1},(\xi_n)_{n\geq 1})$ is the
optimal strategy then we have: $$
\begin{array}{l}
Y^{i,x}_0=E[\int_0^{+\infty}e^{-rs}\psi_{u_s}(X_s^{x})ds-\sum_{k\geq
1} e^{-r\tau_k}g_{u_{\tau_{k-1}}u_{\tau_k}}(X^{x}_{\tau_{k}})].
\end{array}
$$
Taking into account that $g_{ij}\geq \frac{1}{\alpha} >0$ for any
$i\neq j$ we obtain:
$$\begin{array}{ll}
\frac{1}{\alpha}E[\sum_{k=1,n} e^{-r\tau_k}]+ Y^{i,x}_0&\leq
E[\int_0^{+\infty}e^{-rs}\psi_{u_s}(X_s^{x})ds-\sum_{k\geq n+1}
e^{-r\tau_k}g_{u_{\tau_{k-1}}u_{\tau_k}}(X^{x}_{\tau_{k}})].
\end{array}
$$
But for any $k\le n$, $e^{-r\tau_n}\leq e^{-r\tau_k}$
then:$$\begin{array}{ll} \frac{n}{\alpha} E[e^{-r\tau_n}]+
Y^{i,x}_0&\leq
E[\int_0^{+\infty}e^{-rs}\psi_{u_s}(X_s^{x})ds-\sum_{k\geq n+1}
e^{-r\tau_k}g_{u_{\tau_{k-1}}u_{\tau_k}}(X^{x}_{\tau_{k}})]\\{}&\leq
E[\int_0^{+\infty}e^{-rs}\psi_{u_s}(X_s^{x})ds].
\end{array}
$$and then
$$
\begin{array}{ll}
\frac{n}{\alpha} E[e^{-r\tau_n}] &\leq E[\int_0^{+\infty}
e^{-rs}\mid\psi_{u_s}(X_s^{x})\mid ds]-Y^{i,x}_0\\{}&\leq
E[\int_0^{+\infty} e^{-rs}\mid\psi_{u_s}(X_s^{x})\mid
ds]-Y^{0,i,x}_0.
\end {array}
$$
Finally taking into account the facts that $\psi_i$ and $Y^{0,i,x}$
are of polynomial growth, estimate (\ref{estimat1}) for $X^{x}$ and
$\bf H4$ to obtain the desired result. Note that the polynomial
growth of $Y^{0,i,x}$ stems from Proposition 3. $\Box$ \brm The
estimate (\ref{estiopt}) is also valid for the optimal strategy if
at the initial time the state of the plant is an arbitrary $i\in
{\cal I }$. $\Box$ \erm
\medskip


We are now ready to give the main result of this article. \beth The
functions $(v^1,...,v^m):\R^k\rightarrow \R$ are continuous and
solution in viscosity sense of the system of variational
inequalities with inter-connected obstacles (\ref{sysvi1}).\eeth
$Proof$: First let us focus on continuity and let us show that $v1$
is continuous. The same proof will be valid for the continuity of
the other functions $v^i$ ($i=2,...,m$). First the characterization
(\ref{eq27}) implies that:
$$
Y^{1,x}_0=\sup_{(\delta,\xi)\in {\cal
D}}E[\int_0^{+\infty}e^{-rs}\psi_{u_s}(X^{x}_s)ds -\sum_{n\geq 1}
e^{-r\tau_n}g_{u_{\tau_{n-1}}u_{\tau_n}}(X^{x}_{\tau_{n}})]
$$
On the other hand an optimal strategy $(\delta^*,\xi^*)$ exists and
satisfies the estimates (\ref{estiopt}) with the same constant $C$.
Next let $\epsilon >0$ and $x'\in B(x,\epsilon)$ and let us consider
the following set of strategies:$$ \tilde
D:=\{(\delta,\xi)=((\tau_n)_{n\geq 1}, (\xi_n)_{n\geq 0}) \in {\cal
D} \mbox{ such that } \forall n\geq 1, E[e^{-r\tau_n}] \leq
\frac{C(1+(\epsilon +|x|^\gamma))}{n}\}.$$ Therefore the strategy
$(\delta^*,\xi^*)$ belongs to $\tilde D$ and then we have:

$$\begin{array}{ll} Y^{1,x}_0&=\sup_{(\delta,\xi)\in {\tilde
D}}E[\int_0^{+\infty}e^{-rs}\psi_{u_s}(X^{x}_s)ds -\sum_{n\geq
1} e^{-ru_{\tau_n}}g_{u_{\tau_{n-1}}u_{\tau_n}}(X^{x}_{\tau_{n}}) ]\\
{}&=sup_{(\delta,u) \in {\tilde D}}
E[\int_{0}^{\tau_n}e^{-rs}\psi_{u_s}(X^{x}_s)ds
\\{}&\qquad\qquad\qquad-\sum_{1\leq k \leq n} e^{-ru_{\tau_k}}g_{u_{\tau_{k
-1}}u_{\tau_k}}(X^{x}_{{\tau}_{k}})
+e^{-r\tau_n}Y^{u_{\tau_n},x}_{\tau_n}]\end{array}
$$
and
$$\begin{array}{ll} Y^{1,x'}_0&=\sup_{(\delta,\xi)\in {\tilde
D}}E[\int_0^{+\infty}e^{-rs}\psi_{u_s}(X^{x'}_s)ds -\sum_{n\geq
1}e^{-ru_{\tau_n}} g_{u_{\tau_{n-1}}u_{\tau_n}}(X^{x'}_{\tau_{n}})]\\
{}&=sup_{(\delta,u) \in {\tilde D}}
E[\int_{0}^{\tau_n}e^{-rs}\psi_{u_s}(X^{x'}_s)ds
\\{}&\qquad\qquad\qquad-\sum_{1\leq k \leq n}e^{-ru_{\tau_k}} g_{u_{\tau_{k
-1}}u_{\tau_k}}(X^{x'}_{{\tau}_{k}})
+e^{-r\tau_n}Y^{u_{\tau_n},x'}_{\tau_n}]\end{array}
$$

The second equalities it due to the dynamical programming principle.
It follows that:\be \label{eqcont}
\begin{array}{lll}|Y^{1,x'}_0-Y^{1,x}_0|&\leq
sup_{(\delta,u) \in {\tilde D}}
E[\int_{0}^{\tau_n}e^{-rs}|\psi_{u_s}(X^{x'}_s)-\psi_{u_s}(X^{x}_s)|ds\\
{}&\qquad +\sum_{1\leq k \leq n}e^{-ru_{\tau_k}} |g_{u_{\tau_{k
-1}}u_{\tau_k}}(X^{x'}_{\tau_{k}}) -g_{u_{\tau_{k
-1}}u_{\tau_k}}(X^{x}_{\tau_{k}})| \\{}&\qquad+e^{-r\tau_n}|Y^{u_{\tau_n},x'}_{\tau_n}-Y^{u_{\tau_n},x}_{\tau_n}|]\\
{}&\leq
E[\int_{0}^{+\infty}\max_{j=1,m}e^{-rs}|\psi_{j}(X^{x'}_s)-\psi_{j}(X^{x}_s)|ds\\
{}&\qquad +n\max_{i\neq j\in {\cal I}}\{\sup_{s\geq
0}e^{-rs}|g_{ij}(X^{x'}_{s }) -g_{ij}(X^{x}_{s})|\}]
\\{}&\qquad +sup_{(\delta,u) \in {\tilde D}}(E[e^{-2r\tau_n}])^{\frac{1}{2}}(2E[(Y^{u_{\tau_n},x'}_{\tau_n})2+(Y^{u_{\tau_n},x}_{\tau_n})2])^{\frac{1}{2}}.
\end{array}
\ee In the right-hand side of (\ref{eqcont}) the first and the
second term
converges to $0$ as $x'\rightarrow x$. 

Now let us focus on the last one. Since $(\delta,u)\in \tilde D$
then:
$$\begin{array}{ll}sup_{(\delta,u) \in {\tilde D}}(E[e^{-2r\tau_n}])^{\frac{1}{2}}(2E[(Y^{u_{\tau_n},x'}_{\tau_n})2+(Y^{u_{\tau_n},x}_{\tau_n})2]
)^{\frac{1}{2}}&\leq sup_{(\delta,u) \in {\tilde
D}}(E[e^{-r\tau_n}])^{\frac{1}{2}}(2E[(Y^{u_{\tau_n},x'}_{\tau_n})2+(Y^{u_{\tau_n},x}_{\tau_n})2]
)^{\frac{1}{2}}\\&\leq n^{-\frac{1}{2}}\sup_{(\delta,u) \in {\tilde
D}}(2E[(Y^{u_{\tau_n},x'}_{\tau_n})2+(Y^{u_{\tau_n},x}_{\tau_n})2])^{\frac{1}{2}}\\
{}&\leq Cn^{-\frac{1}{2}}(1+|x|^\gamma+|x'|^\gamma)\end{array}$$
where $C$ an appropriate constant which comes from the polynomial
growth of $\psi_i$, $i\in {\cal I}$, estimate (\ref{estimat1}) for
the process $X^{x}$ and inequality (\ref{croi}). Going back now to
(\ref{eqcont}), taking the limit as $x'\rw x$ to obtain:
$$
\lim_{x'\rw x}|Y^{1,x'}_0 - Y^{1,x}_0|\leq
Cn^{-\frac{1}{2}}(1+2|x|^\gamma).$$As $n$ is arbitrary then putting
$n\rw +\infty$ to obtain:
$$
Y^{1,x'}_0 \rightarrow Y^{1,x}_0.$$ Therefore $v^1$ is continuous.
In the same way we can show that $v^2$,...,$v^m$ are continuous. As
they are of polynomial growth then taking into account Theorem 2 to
obtain that $(v^1,\dots,v^m)$ is a viscosity solution for the system
of variational inequalities with inter-connected obstacles
(\ref{sysvi1}). $\Box$
\section{Uniqueness of the solution of the system} We are going now to address the
question of uniqueness of the viscosity solution of the system
(\ref{sysvi1}). We have the following:

\beth \label{uni}The solution in viscosity sense of the system of
variational inequalities with inter-connected obstacles
(\ref{sysvi1}) is unique in the space of continuous functions on
$R^k$ which satisfy a polynomial growth condition, i.e., in the
space
$$\begin{array}{l}{\cal C}:=\{\varphi:  \R^k\rightarrow
\R, \mbox{ continuous and for any }\\\qquad \qquad\qquad x, \,
|\varphi(x)|\leq C(1+|x|^\gamma) \mbox{ for some constants } C\quad
\mbox{and}\quad \gamma\}.\end{array}$$ \eeth {\it Proof}. We will
show by contradiction that if $u_1,...,u_m$ and $w_1,...,w_m$ are a
subsolution and a supersolution respectively for (\ref{sysvi1}) then
for any $i=1,...,m$, $u_i\leq w_i$. Therefore if we have two
solutions of (\ref{sysvi1}) then they are obviously equal. Actually
for some $R>0$ suppose there exists $(x_0,i_0)\in B_R\times {\cal
I}$ $(B_R := \{x\in \R^k; |x|\leq R\})$ such that:
\begin{equation}
\label{comp-equ}
\max\limits_{(x,i)}(u_i(x)-w_i(x))=u_{i_0}(x_0)-w_{i_0}(x_0)=\eta>0.
\end{equation}
 Then, for a small
$\epsilon>0$, and $\theta,\lambda \in(0,1)$ small enough, let us
define:
\begin{equation}
\label{phi}
\Phi^i_{\epsilon}(x,y)=u_{i}(x)-(1-\lambda)w_{i}(y)-\frac{1}{2\epsilon}|x-y|^{2\gamma}
-\theta(|x-x_0|^{2\gamma +2}+|y-x_0|^{2\gamma+2}).
\end{equation}
By the polynomial growth assumption on $u_i$ and $w_i$, there exists
a $(x_{\epsilon},y_{\epsilon},i_\epsilon)\in  B_R \times B_R \times
{\cal I}$, such that:
$$\Phi^{i_\epsilon}_{\epsilon}(x_{\epsilon},y_{\epsilon})=\max\limits_{(x,y,i)}\Phi^i_{\epsilon}(x,y).$$
On the other hand, from
$2\Phi^{i_\epsilon}_{\epsilon}(x_{\epsilon},y_{\epsilon})\geq
\Phi^{i_\epsilon}_{\epsilon}(x_\epsilon,x_\epsilon)+\Phi^{i_\epsilon}_{\epsilon}(y_\epsilon,y_\epsilon)$,
we have
\begin{equation}
\begin{array}{ll}
\frac{1}{2\epsilon}|x_\epsilon -y_\epsilon|^{2\gamma} &\leq (u_{i_\epsilon}(x_\epsilon)-u_{i_\epsilon}(y_\epsilon))+(1-\lambda)(w_{i_\epsilon}(x_\epsilon)-w_{i_\epsilon}(y_\epsilon))\\
&\leq \sum \limits_{i\in{\cal
I}}|u_{i}(x_\epsilon)-u_{i}(y_\epsilon)|+(1-\lambda)\sum
\limits_{i\in{\cal I}}|w_{i}(x_\epsilon)-w_{i}(y_\epsilon)|
\end{array}
\end{equation}

and consequently $\frac{1}{2\epsilon}|x_\epsilon
-y_\epsilon|^{2\gamma}$ is bounded, and as $\epsilon\rightarrow 0$,
$|x_\epsilon -y_\epsilon|\rightarrow 0$. Since $u_{i}$
and $w_{i}$ are uniformly continuous on $ B_R$, then $\frac{1}{2\epsilon}|x_\epsilon -y_\epsilon|^{2\gamma}\rightarrow 0$ as $\epsilon\rightarrow 0.$\\
Since
 $$u_{i_0}(x_0)-(1-\lambda)w_{i_0}(x_0) \leq \Phi^{i_\epsilon}_{\epsilon}(x_{\epsilon},y_{\epsilon})\leq u_{i_\epsilon}(x_\epsilon)-(1-\lambda)w_{i_\epsilon}(y_\epsilon),$$
 it follow as $\lambda \rightarrow 0$ and the continuity of $u_i$ and $w_i$ that, up to a subsequence,
 \begin{equation}\label{subsequence}
 (x_\epsilon,y_\epsilon,i_\epsilon)\rightarrow (x_0,x_0,i_0).
 \end{equation}

 We now claim that:
\begin{equation}
\label{visco-comp1} u_{i_\epsilon}(x_\epsilon)-
\max\limits_{j\in{\cal I}^{-i_\epsilon}}\{-g_{i_{\epsilon}
j}(x_\epsilon)+u_j(x_\epsilon)\} > 0.
\end{equation}
Indeed if
$$u_{i_\epsilon}(x_\epsilon)- \max\limits_{j\in{\cal I}^{-i_\epsilon}}\{-g_{i_{\epsilon} j}(x_\epsilon)+u_j(x_\epsilon)\} \leq
0$$ then there exists $k \in {\cal I}^{-i_\epsilon}$ such that:
$$u_{i_\epsilon}(x_\epsilon) \leq -g_{i_{\epsilon}
k}(x_\epsilon)+u_k(x_\epsilon).$$ From the supersolution property of
$w_{i_\epsilon}(y_\epsilon)$, we have
$$ w_{i_\epsilon}(y_\epsilon)\geq \max\limits_{j\in{\cal I}^{-i_\epsilon}}(-g_{i_{\epsilon} j}(y_\epsilon)+w_j(y_\epsilon)) $$
then
$$ w_{i_\epsilon}(y_\epsilon)\geq -g_{i_{\epsilon} k}(y_\epsilon)+w_k(y_\epsilon).$$
It follows that:
$$u_{i_\epsilon}(x_\epsilon)- (1-\lambda)w_{i_\epsilon}(y_\epsilon) -(u_{k}(x_\epsilon)-(1-\lambda)w_{k}(y_\epsilon))\leq (1-\lambda)g_{i_{\epsilon} k}(y_\epsilon)-g_{i_{\epsilon} k}(x_\epsilon).$$
Now since $g_{ij}\geq \alpha >0$, for every $i\neq j$, and taking
into account of (\ref{phi}) to obtain:
$$
\begin{array}{ll}
\Phi^{i_\epsilon}_{\epsilon}(x_{\epsilon},y_{\epsilon})-\Phi^{k}_{\epsilon}(x_{\epsilon},y_{\epsilon})&<
-\alpha \lambda
+g_{i_{\epsilon} k}(y_\epsilon)-g_{i_{\epsilon} k}(x_\epsilon)\\
\end{array}
$$
But this contradicts the definition of $i_\epsilon$, since
$g_{i_{\epsilon} k}$ is uniformly continuous on $B_R$ and the claim
(\ref{visco-comp1}) holds.

Next let us denote
\begin{equation}
\varphi_{\epsilon}(x,y)=\frac{1}{2\epsilon}|x-y|^{2\gamma}
+\theta(|x-x_0|^{2\gamma +2}+|y-x_0|^{2\gamma+2}).
\end{equation}
Then we have: \be \left\{
\begin{array}{lll}\label{derive}
D_{x}\varphi_{\epsilon}(t,x,y)=
\frac{\gamma}{\epsilon}(x-y)|x-y|^{2\gamma-2} +\theta(2\gamma + 2)
(x-x_0)|x-x_0|^{2\gamma}, \\
D_{y}\varphi_{\epsilon}(t,x,y)=
-\frac{\gamma}{\epsilon}(x-y)|x-y|^{2\gamma-2} +
\theta(2\gamma + 2)(y-y_0)|y-y_0|^{2\gamma},\\\\
B(t,x,y)=D_{x,y}^{2}\varphi_{\epsilon}(t,x,y)=\frac{1}{\epsilon}
\begin{pmatrix}
a_1(x,y)&-a_1(x,y) \\
-a_1(x,y)&a_1(x,y)
\end{pmatrix}+ \begin{pmatrix}
a_2(x)&0 \\
0&a_2(y)
\end{pmatrix} \\\\
\mbox{ with } a_1(x,y)=\gamma|x-y|^{2\gamma-2}I+\gamma(2\gamma -2)(x-y)(x-y)^* |x-y|^{2\gamma-4} \mbox{ and }\\
a_2(x)=\theta(2\gamma + 2)|x-x_0|^{2\gamma}I+2\theta \gamma(2\gamma
+ 2)(x-x_0)(x-x_0)^* |x-x_0|^{2\gamma-2 }.
\end{array}
\right. \ee Taking into account (\ref{visco-comp1}) then applying
the result by Crandall et al. (Theorem 3.2, {\cite{[CIL]}) to the
function $$ u_{i}(x)-(1-\lambda)w_{i}(y)-\varphi_{\epsilon}(x,y) $$
at the point $(x_\epsilon,y_\epsilon)$, for any $\epsilon_1 >0$, we
can find
 $X,Y \in S_k$, such that:

\be \label{lemmeishii} \left\{
\begin{array}{lllll}
(\frac{\gamma}{\epsilon}(x_\epsilon-y_\epsilon)|x_\epsilon-y_\epsilon|^{2\gamma-2}
+\theta(2\gamma + 2) (x_\epsilon-x_0)|x_\epsilon-x_0|^{2\gamma},X)
\in J^{2,+}(u_{i_\epsilon}(x_\epsilon)),\\
(\frac{\gamma}{\epsilon}(x_\epsilon-y_\epsilon)|x_\epsilon-y_\epsilon|^{2\gamma-2}
- \theta(2\gamma +
2)(y_\epsilon-y_0)|y_\epsilon-y_0|^{2\gamma},Y)\in J^{2,-}
((1-\lambda)w_{i_\epsilon}(y_\epsilon)),\\
-(\frac{1}{\epsilon_1}+||B(x_\epsilon,y_\epsilon)||)I\leq
\begin{pmatrix}
X&0 \\
0&-Y
\end{pmatrix}\leq B(x_\epsilon,y_\epsilon)+\epsilon_1 B(x_\epsilon,y_\epsilon)2.
\end{array}
\right. \ee Taking now into account (\ref{visco-comp1}), and the
definition of viscosity solution, we get:
$$\begin{array}{l}ru_{i_\epsilon}(x_\epsilon)-\frac{1}{2}Tr[\sigma^{*}(x_\epsilon)X\sigma(x_\epsilon)]-\langle\frac{\gamma}{\epsilon}
(x_\epsilon-y_\epsilon)|x_\epsilon-y_\epsilon|^{2\gamma-2}
\\\qquad\qquad\qquad\qquad\qquad +\theta(2\gamma + 2)
(x_\epsilon-x_0)|x_\epsilon-x_0|^{2\gamma},b(x_\epsilon)\rangle-\psi_{i_\epsilon}(x_\epsilon)\leq
0 \mbox{ and
}\\r(1-\lambda)w_{i_\epsilon}(y_\epsilon)-\frac{1}{2}Tr[\sigma^{*}(y_\epsilon)Y\sigma(y_\epsilon)]-\langle
\frac{\gamma}{\epsilon}
(x_\epsilon-y_\epsilon)|x_\epsilon-y_\epsilon|^{2\gamma-2}
\\\qquad\qquad\qquad\qquad\qquad -\theta(2\gamma + 2)
(y_\epsilon-x_0)|y_\epsilon-x_0|^{2\gamma},b(y_\epsilon)\rangle-(1-\lambda)\psi_{i_\epsilon}(y_\epsilon)\geq
0\end{array}$$ which implies that:
\begin{equation}
\begin{array}{llll}
\label{viscder}
&ru_{i_\epsilon}(x_\epsilon)-r(1-\lambda)w_{i_\epsilon}(y_\epsilon)\leq \frac{1}{2}Tr[\sigma^{*}(x_\epsilon)X\sigma(x_\epsilon)-\sigma^{*}(y_\epsilon)Y\sigma(y_\epsilon)]\\
&\qquad + \langle\frac{\gamma}{\epsilon}
(x_\epsilon-y_\epsilon)|x_\epsilon-y_\epsilon|^{2\gamma-2},b(x_\epsilon)-b(y_\epsilon)\rangle\\&\qquad+\langle
\theta(2\gamma + 2)
(x_\epsilon-x_0)|x_\epsilon-x_0|^{2\gamma},b(x_\epsilon)\rangle
+\langle \theta(2\gamma + 2)
(y_\epsilon-x_0)|y_\epsilon-x_0|^{2\gamma},b(y_\epsilon)\rangle
\\&\qquad+\psi_{i_\epsilon}(x_\epsilon)-(1-\lambda)\psi_{i_\epsilon}(y_\epsilon).
\end{array}
\end{equation}
But from (\ref{derive}) there exist two constants $C$ and $C_1$ such
that:
$$||a_1(x_\epsilon,y_\epsilon)||\leq C|x_\epsilon - y_\epsilon|^{2\gamma -2} \mbox{ and }
(||a_2(x_\epsilon)||\vee ||a_2(y_\epsilon)||)\leq C_1\theta .$$ As
$$B= B(x_\epsilon,y_\epsilon)= \frac{1}{\epsilon}
\begin{pmatrix}
a_1(x_\epsilon,y_\epsilon)&-a_1(x_\epsilon,y_\epsilon) \\
-a_1(x_\epsilon,y_\epsilon)&a_1(x_\epsilon,y_\epsilon)
\end{pmatrix}+ \begin{pmatrix}
a_2(x_\epsilon)&0 \\
0&a_2(y_\epsilon)
\end{pmatrix}$$
then
$$B\leq \frac{1}{\epsilon}
\begin{pmatrix}
I&-I \\
-I&I
\end{pmatrix}+ C_1\theta I.$$
It follows that:
\begin{equation}
B+\epsilon_1 B2 \leq C(\frac{1}{\epsilon}|x_\epsilon -
y_\epsilon|^{2\gamma -2}+ \frac{\epsilon_1}{\epsilon2}|x_\epsilon -
y_\epsilon|^{4\gamma -4})\begin{pmatrix}
I&-I \\
-I&I
\end{pmatrix}+ C_1\theta I
\end{equation}
where $C$ and $C_1$ which hereafter may change from line to line.
Choosing now $\epsilon_1=\epsilon$, yields the relation
\begin{equation}
\label{ineg_matreciel} B+\epsilon_1 B2 \leq
\frac{C}{\epsilon}(|x_\epsilon - y_\epsilon|^{2\gamma
-2}+|x_\epsilon - y_\epsilon|^{4\gamma -4})\begin{pmatrix}
I&-I \\
-I&I
\end{pmatrix}+ C_1\theta I.
\end{equation}
Now, from $\bf H1$, (\ref{lemmeishii}) and (\ref{ineg_matreciel}) we
get:
$$\frac{1}{2}Tr[\sigma^{*}(x_\epsilon)X\sigma(x_\epsilon)-\sigma^{*}(y_\epsilon)
Y\sigma(y_\epsilon)]\leq \frac{C}{\epsilon}(|x_\epsilon -
y_\epsilon|^{2\gamma}+|x_\epsilon - y_\epsilon|^{4\gamma -2}) +C_1
\theta(1+|x_\epsilon|^2+|y_\epsilon|^2).$$ Next $$
\langle\frac{\gamma}{\epsilon}
(x_\epsilon-y_\epsilon)|x_\epsilon-y_\epsilon|^{2\gamma-2},b(x_\epsilon)-b(y_\epsilon)\rangle
\leq \frac{C2}{\epsilon}|x_\epsilon - y_\epsilon|^{2\gamma}$$ and
finally,
$$\langle
\theta(2\gamma + 2)
(x_\epsilon-x_0)|x_\epsilon-x_0|^{2\gamma},b(x_\epsilon)\rangle \leq
\theta C(1+|x_\epsilon|)|x_\epsilon -x_0|^{2\gamma+1}
$$
$$
\langle \theta(2\gamma + 2)
(y_\epsilon-x_0)|y_\epsilon-x_0|^{2\gamma},b(y_\epsilon)\rangle \leq
\theta C(1+|y_\epsilon|)|y_\epsilon -x_0|^{2\gamma+1}.$$ So that by
plugging into (\ref{viscder}) we obtain:
$$\begin{array}{l}ru_{i_\epsilon}(x_\epsilon)-r(1-\lambda)w_{i_\epsilon}(y_\epsilon) \leq \frac{C}{\epsilon}(|x_\epsilon -
y_\epsilon|^{2\gamma}+|x_\epsilon - y_\epsilon|^{4\gamma -2}) +C_1
\theta(1+|x_\epsilon|^2+|y_\epsilon|^2)+\frac{C2}{\epsilon}|x_\epsilon
- y_\epsilon|^{2\gamma}+\\\qquad \qquad \theta
C(1+|x_\epsilon|)|x_\epsilon -x_0|^{2\gamma+1}+\theta
C(1+|y_\epsilon|)|y_\epsilon
-x_0|^{2\gamma+1}+\psi_{i_\epsilon}(x_\epsilon)-(1-\lambda)\psi_{i_\epsilon}(y_\epsilon).\end{array}$$
By sending $\epsilon\rightarrow0$, $\lambda\rightarrow0$, $\theta
\rightarrow0$ and taking into account of the continuity of
$\psi_{i_\epsilon}$, we obtain $u_{i_0}(x_0)-w_{i_0}(x_0)<0$ which
is a contradiction. The proof of Theorem \ref{uni} is now complete.
$\Box$
\medskip

As a by-product we have the following Corollary: \bcor Let
$(v^1,...,v^m)$ be a viscosity solution of (\ref{sysvi1}) which
satisfies a polynomial growth condition then for $i=1,...,m$ and
$(t,x)\in \R^k$, $$ v^i(x)= \sup_{(\delta,\xi)\in {\cal
D}^i_0}E[\integ{0}{+\infty}e^{-rs}\psi_{u_s}(X^{x}_s)ds -\sum_{n\geq
1} e^{-r\tau_n}g_{u_{\tau_{n-1}}u_{\tau_n}}(X^{x}_{\tau_{n}}) ].
$$\\
\section {Numerical results} We consider now some numerical examples
of the optimal switching problem (\ref{sysvi1}).\\
{\bf Example1}: In this example we consider an optimal switching
problem with two
modes, where\\
$r=100$, $b=x$, $\sigma=\sqrt{2}x$,
$g_{12}(x)={\frac{1}{2}}|x|+0.1$,
$g_{21}(t,x)=|x|+0.48$, $\psi_1(x)={\frac{1}{2}}x^2-0.3x+1$, $\psi_2(t,x)=x^2+1$.\\
\begin{figure}[h]
    \begin{center}
      \includegraphics{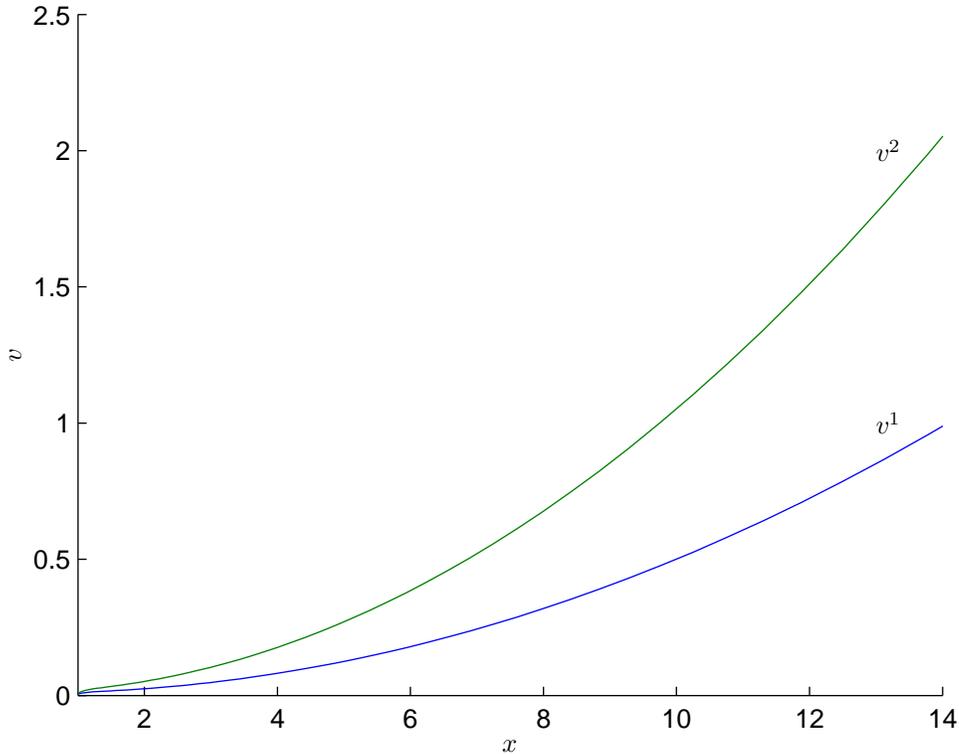}
    \end{center}
    \caption{Curves of $v^2$ and $v^1$.}
\end{figure}\\
{\bf Example2}: We now consider the case of 3 modes where $r=100$,
$b=x$, $\sigma=\sqrt{2}x$, $g_{12}(t,x)=0.5|x|+1$,
$g_{13}(t,x)=x^2+0.5$, $g_{21}(t,x)=|x|+4$, $g_{23}(t,x)=|x|+5$,
$g_{31}(t,x)=0.001|x|+0.1$, $g_{32}(t,x)=x^2+|x|+0.5$,
$\psi_1(t,x)=x+1$, $\psi_2(t,x)=-x-2$ and finally
$\psi_3(t,x)=-x-2$.
\begin{figure}[h]
    \begin{center}
      \includegraphics{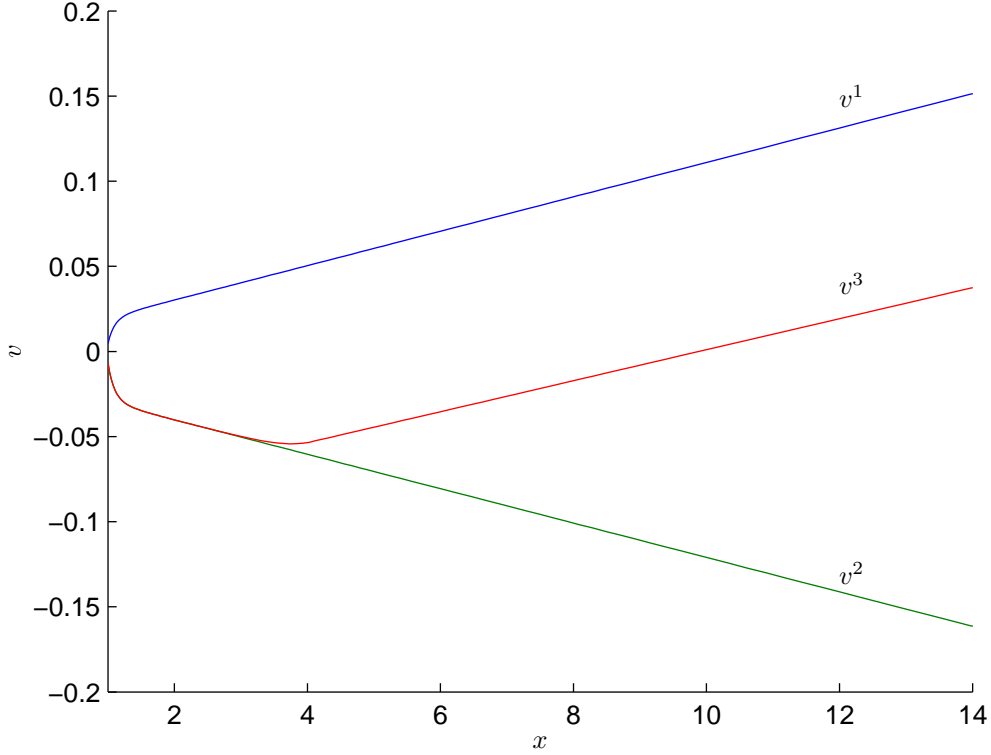}
    \end{center}
    \caption{Curves of $v^1$, $v^3$ and $v^2$.}
\end{figure}

 \ecor

 \no {\bf Acknowledgement}: The author thanks gratefully Prof. S.
Hamad\`ene for the fructuous discussions during the preparation of
this paper.$\Box$
\section*{ Appendix: proof of Theorem 1}

The proof consists in showing that for any $t\leq T,$ $Y^i_t$, as
defined by (\ref{eqvt}), is nothing but the expected total profit or
the value function of the optimal problem, given that the system is
in mode $i$ at time $t$. More precisely,
$$ e^{-rt}Y^i_t=\esssup_{(\delta,u)\in {\cal D}_t}E[\int_t^{+\infty} e^{-rs}\psi_i(X_s)ds
- \sum_{ k\geq1} e^{-r\tau_k}g_{u_{\tau_{k
-1}}u_{\tau_k}}(X_{{\tau_{k}}})|\cF_t],$$

where ${\cal D}_t$ is the set of strategies such that $\tau_1\geq
t$, P-a.s. if at time $t$ the system is in the mode i.

Let us admit for a moment the following Lemma.
\begin{lem}\label{env-sne} For every $t\geq \tau^*_1$.

\begin{equation}\label{env1} e^{-rt}
Y^{u_{\tau^*_1}}_t=\esssup_{\tau \geq t
}E[\int_t^{\tau}e^{-rs}\psi_{u_{\tau^*_1}}(X_s)ds
+e^{-r\tau}\max\limits_{j\in {\cal
I}^{-u_{\tau^*_1}}}(-g_{ij}(X_{\tau})+Y^j_{\tau})|\cF_t].\Box
\end{equation}

\end{lem} From properties of the Snell envelope and at time $t=0$ the system is in mode $1$, we have:

$$
\begin{array}{ll}
 Y^1_0&=E[\int_0^{\tau^*_1}e^{-rs}\psi_1(X_s)ds
+e^{-r\tau^*_1}\max\limits_{j\in {\cal
I}^{-i}}(-g_{ij}(X_{\tau^*_1})+Y^j_{\tau^*_1})]\\
&=E[\int_0^{\tau^*_1}e^{-rs}\psi_1(X_s)ds
+e^{-r\tau^*_1}(-g_{iu_{\tau^*_1}}(X_{\tau^*_1})+Y^{u_{\tau^*_1}}_{\tau^*_1})].
\end{array}$$
Now from Lemma 2 and the definition of $\tau^*_2$ we have:
$$
\begin{array}{ll}
e^{-r\tau^*_1}Y^{u_{\tau^*_1}}_{\tau^*_1}&=E[\int_{\tau^*_1}^{\tau^*_2}e^{-rs}\psi_{u_{\tau^*_1}}(X_s)ds
+e^{-r\tau^*_2}\max\limits_{j\in {\cal
I}^{-u_{\tau^*_1}}}(-g_{u_{\tau^*_1}j}(X_{\tau^*_2})+Y^j_{\tau^*_2})|\cF_{\tau^*_1}]\\
&=E[\int_{\tau^*_1}^{\tau^*_2}e^{-rs}\psi_{u_{\tau^*_1}}(X_s)ds
+e^{-r\tau^*_2}(-g_{u_{\tau^*_1}u_{\tau^*_2}}(X_{\tau^*_2})+Y^{u_{\tau^*_2}}_{\tau^*_2})|\cF_{\tau^*_1}].
\end{array}
$$

It implies that
$$
\begin{array}{lll}
Y^1_0&=E[\int_0^{\tau^*_1}e^{-rs}\psi_1(X_s)ds
-e^{-r\tau^*_1}g_{iu_{\tau^*_1}}(X_{\tau^*_1})\\
&+ E[\int_{\tau^*_1}^{\tau^*_2}e^{-rs}\psi_{u_{\tau^*_1}}(X_s)ds
+e^{-r\tau^*_2}(-g_{u_{\tau^*_1}u_{\tau^*_2}}(X_{\tau^*_2})+Y^{u_{\tau^*_2}}_{\tau^*_2})|\cF_{\tau^*_1}]]\\
&=E[\int_0^{\tau^*_1}e^{-rs}\psi_1(X_s)ds+
\int_{\tau^*_1}^{\tau^*_2}e^{-rs}\psi_{u_{\tau^*_1}}(X_s)ds
-e^{-r\tau^*_1}g_{iu_{\tau^*_1}}(X_{\tau^*_1})-e^{-r\tau^*_2}g_{u_{\tau^*_1}u_{\tau^*_2}}(X_{\tau^*_2})
\\&+ e^{-r\tau^*_2}Y^{u_{\tau^*_2}}_{\tau^*_2}].
\end{array}
$$
 Therefore
$$Y^1_0=E[\int_0^{\tau^*_2}e^{-rs}\psi(X_s,u_s)ds -e^{-r\tau^*_1}g_{iu_{\tau^*_1}}(X_{\tau^*_1})-e^{-r\tau^*_2}g_{u_{\tau^*_1}u_{\tau^*_2}}(X_{\tau^*_2})
 + e^{-r\tau^*_2}Y^{u_{\tau^*_2}}_{\tau^*_2}],$$ since between 0 and
$\tau^*_1$ (resp. $\tau^*_1$ and $\tau^*_2$) the production is in
regime $1$ (resp. regime $u_{\tau^*_1}$) and then $u_t=1$ (resp.
$u_t=u_{\tau^*_1}$) which implies that
$$\integ{0}{\tau^*_2}e^{-rs}\psi(X_s,u_s)ds=\integ{0}{\tau^*_1}
e^{-rs}\psi_1(X_s)ds+\integ{\tau^*_1}{\tau^*_2}
e^{-rs}\psi_{u_{\tau^*_1}}(X_s)ds.$$ Now repeating this reasoning as
many times as necessary we obtain that for any $n\geq 0,$
$$\begin{array}{l}
Y^1_0= E[\integ{0}{\tau^*_{n}}e^{-rs}\psi(X_s,u_s)ds -\sum_{1\leq
k\leq n} e^{-r\tau^*_k}g_{u_{\tau^*_{k
-1}}u_{\tau^*_k}}(X_{{\tau^*_k}})+
e^{-r\tau^*_n}Y^{u_{\tau^*_n}}_{\tau^*_{n}}].
\end{array}$$

Then, the strategy $(\delta^*,u^*)$ verify $E[\sum_{n\geq0}
e^{-r\tau^*_{n}}]<+\infty$, otherwise $Y^{1}_0$ would be equal to
$-\infty$ contradicting the assumption that the processes $Y^i$
belong to ${\cal S}^2$. Therefore, taking the limit as $n
\rightarrow
+\infty$ we obtain $Y^1_0 =J(\delta^*,u^*)$.\\
To complete the proof it remains to show that the strategy
$(\delta^*,u^*)$ it is optimal i.e. $J(\delta^*,u^*)\geq
J(\delta,u)$ for any $(\delta,u)\in \cal D$.\\
The definition of the Snell envelope yields
$$
\begin{array}{ll}
 Y^1_0&\geq E[\int_0^{\tau_1}e^{-rs}\psi_1(X_s)ds +e^{-r\tau_1}\max\limits_{j\in {\cal
I}^{-1}}(-g_{1j}(X_{\tau_1})+Y^j_{\tau_1})]\\
&\geq E[\int_0^{\tau_1}e^{-rs}\psi_1(X_s)ds
+e^{-r\tau_1}(-g_{1u_{\tau^*_1}}(X_{\tau_1})+Y^{u_{\tau_1}}_{\tau_1})].
\end{array}$$
 But, once more using a similar characterization as (\ref{env1}), we get
$$
\begin{array}{ll}
e^{-r\tau_1}Y^{u_{\tau_1}}_{\tau_1}&\geq
E[\int_{\tau_1}^{\tau_2}e^{-rs}\psi_{u_{\tau_1}}(X_s)ds
+e^{-r\tau_2}\max\limits_{j\in {\cal
I}^{-u_{\tau_1}}}(-g_{u_{\tau_1}j}(X_{\tau_2})+Y^j_{\tau_2})|\cF_{\tau_1}]\\
&\geq E[\int_{\tau_1}^{\tau_2}e^{-rs}\psi_{u_{\tau_1}}(X_s)ds
+e^{-r\tau_2}(-g_{u_{\tau_1}u_{\tau_2}}(X_{\tau_2})+Y^{u_{\tau_2}}_{\tau_2})|\cF_{\tau_1}].
\end{array}
$$
Therefore,
$$
\begin{array}{lll}
Y^1_0&\geq E[\int_0^{\tau_1}e^{-rs}\psi_1(X_s)ds
-e^{-r\tau_1}g_{1u_{\tau_1}}(X_{\tau_1})]\\
&+ E[\int_{\tau_1}^{\tau_2}e^{-rs}\psi_{u_{\tau_1}}(X_s)ds
+e^{-r\tau_2}(-g_{u_{\tau_1}u_{\tau_2}}(X_{\tau_2})+Y^{u_{\tau_2}}_{\tau_2})]\\
&=E[\int_0^{\tau_2}e^{-rs}\psi(X_s,u_s)ds
-e^{-r\tau_1}g_{1u_{\tau_1}}(X_{\tau_1})-e^{-r\tau_2}g_{u_{\tau_1}u_{\tau_2}}(X_{\tau_2})+
e^{-r\tau_2}Y^{u_{\tau_2}}_{\tau_2}].
\end{array}
$$
Repeat this argument $n$ times to obtain
$$\begin{array}{l}
Y^1_0\geq E[\integ{0}{\tau_{n}}e^{-rs}\psi(X_s,u_s)ds -\sum_{1\leq
k\leq n} e^{-r\tau_k}g_{u_{\tau_{k -1}}u_{\tau_k}}(X_{{\tau_k}})+
e^{-r\tau_n}Y^{u_{\tau_n}}_{\tau_{n}}].
\end{array}$$
Finally, taking the limit as $n\rightarrow +\infty$ yields
$$\begin{array}{l}
Y1_0\geq E[\integ{0}{+\infty}e^{-rs}\psi(X_s,u_s)ds -\sum_{ k\geq1}
e^{-r\tau_k}g_{u_{\tau_{k -1}}u_{\tau_k}}(X_{{\tau_k}})].
\end{array}$$
 Hence, the
strategy $(\delta^*,u^*)$ is optimal. We proceed to the proof of Lemma 2.\\
\indent $Proof$ of Lemma 2. From (\ref{eqvt}) we have for any $i\in
{\cal I}$ and $t\geq 0$\be
\begin{array}{l}
e^{-rt}Y^i_t=\esssup_{\tau \geq t}E[\int_t^\tau e^{-rs}\psi_i(X_s)ds
+e^{-r\tau}\max\limits_{j\in {\cal
I}^{-i}}(-g_{ij}(X_\tau)+Y^j_\tau)|\cF_t].
\end{array}
\ee This also means that the process $(e^{-rt}Y^i_t+\int_0^t
e^{-rs}\psi_i(X_s)ds)_{t\geq 0}$ is a supermartingale which
dominates $$(\int_0^te^{-rs}\psi_i(X_s)ds+e^{-rt}\max\limits_{j\in
{\cal I}^{-i}}(-g_{ij}(X_t)+Y^j_t))_{t\geq 0}.$$ This implies that
the process
$(\ind_{[u_{\tau^*_1}=i]}(e^{-rt}Y^i_t+\int_{\tau^*_1}^te^{-rs}\psi_i(X_s)ds))_{t\geq
\tau^*_1}$  is a supermartingale which dominates
$$(\ind_{[u_{\tau^*_1}=i]}(\int_{\tau^*_1}^te^{-rs}\psi_i(X_s)ds+e^{-rt}\max\limits_{j\in {\cal
I}^{-i}}(-g_{ij}(X_t)+Y^j_t))_{t\geq \tau^*_1}.$$ Since ${\cal I}$
is finite, the process $(\sum_{i\in {\cal
I}}\ind_{[u_{\tau^*_1}=i]}(e^{-rt}Y^i_t+\int_{\tau^*_1}^te^{-rs}\psi_i(X_s)ds))_{t\geq
\tau^*_1}$  is also a supermartingale which dominates
$(\sum_{i\in{\cal
I}}\ind_{[u_{\tau^*_1}=i]}(\int_{\tau^*_1}^te^{-rs}\psi_i(X_s)ds+e^{-rt}\max\limits_{j\in
{\cal I}^{-i}}(-g_{ij}(X_t)+Y^j_t))_{t\geq \tau^*_1}.$\\ Thus, the
process $(e^{-rt}Y^{u_{\tau^*_1}}_t+\int_{\tau^*_1}^t
e^{-rs}\psi_{u_{\tau^*_1}}(X_s)ds)_{t\geq \tau^*_1}$ is a
supermartingale which is greater than
$$(\int_{\tau^*_1}^te^{-rs}\psi_{u_{\tau^*_1}}(X_s)ds+e^{-rt}\max\limits_{j\in
{\cal I}^{-u_{\tau^*_1}}}(-g_{u_{\tau^*_1}j}(X_t)+Y^j_t))_{t\geq
\tau^*_1}.$$ To complete the proof it remains to show that it is the
smallest one which has this property and use the characterization of
the Snell envelope see e.g. \cite{[CK], Elka, ham}.\\ Indeed, let
$(Z_t)_{t\geq \tau^*_1}$ be a supermartingale of class $[D]$ such
that, for any $t\geq \tau^*_1$,
$$Z_t \geq \int_{\tau^*_1}^te^{-rs}\psi_{u_{\tau^*_1}}(X_s)ds+e^{-rt}\max\limits_{j\in
{\cal I}^{-u_{\tau^*_1}}}(-g_{u_{\tau^*_1}j}(X_t)+Y^j_t).$$ It
follows that for every $t\geq \tau^*_1$,
$$\ind_{[u_{\tau^*_1}=i]}Z_t \geq \ind_{[u_{\tau^*_1}=i]}(\int_{\tau^*_1}^te^{-rs}\psi_{i}(X_s)ds+e^{-rt}\max\limits_{j\in
{\cal I}^{-i}}(-g_{ij}(X_t)+Y^j_t)).$$ But, the process
$(\ind_{[u_{\tau^*_1}=i]}Z_t)_{t\geq \tau^*_1}$ is a supermartingale
and for every $t\geq \tau^*_1$,
$$\ind_{[u_{\tau^*_1}=i]}e^{-rt}
Y^{i}_t=\esssup_{\tau \geq t
}E[\ind_{[u_{\tau^*_1}=i]}(\int_t^{\tau}e^{-rs}\psi_{i}(X_s)ds
+e^{-r\tau}\max\limits_{j\in {\cal
I}^{-i}}(-g_{ij}(X_{\tau})+Y^j_{\tau}))|\cF_t].$$ It follows that,
for every $t\geq \tau^*_1$, $$\ind_{[u_{\tau^*_1}=i]}Z_t\geq
\ind_{[u_{\tau^*_1}=i]}(e^{-rt} Y^{i}_t +
\int_{\tau^*_1}^te^{-rs}\psi_{i}(X_s)ds).$$ Summing over $i$, we
get, for every  $t\geq \tau^*_1$,
$$Z_t\geq
e^{-rt} Y^{u_{\tau^*_1}}_t +
\int_{\tau^*_1}^te^{-rs}\psi_{u_{\tau^*_1}}(X_s)ds.$$ Hence, the
process $(e^{-rt} Y^{u_{\tau^*_1}}_t +
\int_{\tau^*_1}^te^{-rs}\psi_{u_{\tau^*_1}}(X_s)ds)_{t\geq
\tau^*_1}$ is the Snell envelope of $$(\int_{\tau^*_1}^t
e^{-rs}\psi_{u_{\tau^*_1}}(X_s)ds +e^{-rt}\max\limits_{j\in {\cal
I}^{-u_{\tau^*_1}}}(-g_{u_{\tau^*_1}j}(X_t)+Y^j_{t}))_{t\geq
\tau^*_1},$$  whence Lemma 2.$\Box$



\begin{thebibliography}{99}
\small
\renewcommand{\baselinestretch}{0.3}
\bibitem{[BE]} Bayraktar, E. and Egami, M. (2007): On the One-Dimensional Optimal Switching Problem. {\it Preprint}.

\bibitem{[BO1]} Brekke, K. A. and \O ksendal, B. (1994): Optimal
switching in an economic activity under uncertainty. {\it SIAM J.
Control Optim.} (32), pp. 1021-1036.


\bibitem{[BS]} Brennan, M. J. and Schwartz, E. S. (1985): Evaluating
natural resource investments. {\it J.Business} 58,  pp. 135-137.

\bibitem{[CL]} Carmona, R. and Ludkovski, M. (2005): Optimal
Switching with Applications to Energy Tolling Agreements. {\it
Preprint}.

\bibitem{[CH]} Chen, Z. (1998): Existence and uniqueness for BSDE's
with stopping time, Chinese Science Bulletin, 43, p.96-99.

\bibitem{[CIL]} Crandall, M., Ishii, H. and P.L. Lions (1992) : User's guide to viscosity solutions of
second order partial differential equations, Bull. Amer. Math. Soc.,
27, 1-67.

\bibitem{[CK]} Cvitanic, J. and Karatzas, I (1996): Backward SDEs with
reflection and Dynkin games. {\it Annals of Probability} 24 (4), pp.
2024-2056.

\bibitem{[DEL]} \textsc{Dellacherie, C.} and \textsc{Meyer, P.A.} (1980). Probabilit\'es et Potentiel V-VIII, {\it Hermann,
Paris.}



\bibitem{[DP]} Dixit, A. and  Pindyck, R. S. (1994): Investment under
uncertainty. {\it Princeton Univ. Press.}

\bibitem{[DH]} Djehiche, B. and Hamad\`ene, S (2009): On a finite horizon
Starting and Stopping Problem with Default risk. {\it  to appear in
the International J. of Theoretical and Applied Finance (IJTAF).}

\bibitem{[DHP]} Djehiche, B. Hamad\`ene, S. and Popier, A. (2007): A finite horizon optimal multiple switching problem.
{\it Preprint, Universit\'e du Maine, F.}



\bibitem{[DZ2]} Duckworth, K. and Zervos, M. (2001): A model for
investment decisions with switching costs. {\it Annals of Applied
probability} 11 (1), pp. 239-260.

\bibitem{[EH]} El Asri, B. and Hamad\`ene, S. (2009): The Finite Horizon Optimal Multi-Modes Switching Problem: the Viscosity Solution
Approach, {\it  Applied Mathematics and Optimization, DOI
10.1007/s00245-009-9071-3.}


\bibitem{Elka} El Karoui, N. (1980): Les aspects probabilistes du
contr\^ ole stochastique. {\it Ecole d'\'et\'e de probabilit\'es de
Saint-Flour, Lect. Notes in Math. No 876, Springer Verlag.}


\bibitem{[EKal]} El-Karoui, N. Kapoudjian, C. Pardoux, E. Peng, S.
and Quenez, M. C. (1997): Reflected solutions of backward SDEs and
related obstacle problems for PDEs. {\it Annals of Probability} 25
(2), pp. 702-737.

\bibitem{ham} Hamad\`ene, S. (2002): Reflected BSDEs with discontinuous
barriers. {\it Stochastics and Stochastic Reports 74 (3-4), pp.
571-596.}


\bibitem{[HJ]} Hamad\`ene, S. and  Jeanblanc, M (2007):
On the Starting and Stopping Problem: Application in reversible
investments, {\it Math. of Operation Research, vol.32, No.1,
pp.182-192}.

\bibitem{hib} Hamad\`ene, S. and Hdhiri, I. (2006): On the starting and
stopping problem in the model with jumps. {\it Preprint ,
Universit\'e du Maine, Le Mans, F.}

\bibitem{[HLW]} Hamad\`ene, S. Lepeltier, J-P and Wu, Z. (1999): nfinite Horizon Reflected BSDE's and Applications in Mixed Control and Game Problems. {\it Probability and Mathematical Statistics International Journal vol.19, pp.211-234}



\bibitem{[LPZ]} Ly Vath, V. Pham, H and Zhou, X. (2007): Optimal switching over multiple regimes.
{\it Preprint.}.


\bibitem{[PA]} Pardoux, E. (1999):
Weak convergence and homogenization of semilinear PDEs. {\it Nonlin.
Anal, Dif. Equa. and Control, pp. 503-549.}


%

\bibitem{[RY]} Revuz, D and Yor, M. (1991): Continuous Martingales and
Brownian Motion. {\it Springer Verlag, Berlin.}


\bibitem{[TY]} Tang, S. and Yong, J. (1993): Finite horizon stochastic optimal switching and impulse
controls with a viscosity solution approach.  Stoch. and Stoch.
Reports, 45, 145-176.

%

\bibitem{dz} Zervos, M. (2003): A Problem of Sequential Enty and Exit
Decisions Combined with Discretionary Stopping. {\it SIAM J. Control
Optim.} 42 (2), pp. 397-421.

\end{thebibliography}
\end{document}